\documentclass[11pt]{article}

\usepackage{hyperref}
\usepackage{amsmath}
\usepackage{amssymb}
\usepackage{amsthm}
\usepackage{latexsym}
\usepackage{color}
\usepackage{graphicx}
\usepackage{appendix}
\usepackage{color} 
\usepackage{enumerate}

\usepackage[T1]{fontenc}
\usepackage[english]{babel}
\usepackage{geometry}
\DeclareSymbolFont{calletters}{OMS}{cmsy}{m}{n}
\DeclareSymbolFontAlphabet{\mathcal}{calletters}

\def\be{\begin{eqnarray}}
\def\ee{\end{eqnarray}}

\def\b*{\begin{eqnarray*}}
\def\e*{\end{eqnarray*}}

\newtheorem{Theorem}{Theorem}[part]
\newtheorem{Definition}{Definition}[part]
\newtheorem{Proposition}{Proposition}[part]

\newtheorem{Assumption}{Assumption}[part]
\newtheorem{Lemma}{Lemma}[part]
\newtheorem{Corollary}{Corollary}[part]
\newtheorem{Remark}{Remark}[part]

\makeatletter \@addtoreset{equation}{section}

\@addtoreset{Definition}{section}

\@addtoreset{Theorem}{section}

\@addtoreset{Proposition}{section}

\@addtoreset{Property}{section}

\@addtoreset{Assumption}{section}

\@addtoreset{Corollary}{section}

\@addtoreset{Lemma}{section}

\@addtoreset{Remark}{section}

\@addtoreset{Example}{section}

\newcommand{\No}[1]{\left\|#1\right\|}     
\newcommand{\abs}[1]{\left|#1\right|}     

\def \E{\mathbb{E}}

\def \P{\mathbb{P}}

\def \R{\mathbb{R}}

\def\Ac{{\cal A}}
\def\Bc{{\cal B}}

\def\Ec{{\cal E}}

\def\Uc{{\cal U}}

\def\esup{{\rm ess \, sup}}

\def\x{\times}

\def\={\;=\;}
\def\.{\;.}

\def\eps{\varepsilon}

\def\reff#1{{\rm(\ref{#1})}}

\def\1{{\bf 1}}
\def \ep{\hbox{ }\hfill{ ${\cal t}$~\hspace{-5.1mm}~${\cal u}$   } }

\def \proof{{\noindent \bf Proof. }}
\def \ep{\hbox{ }\hfill$\Box$}

 \def\normeL2#1{\left\|{#1}\right\|_{L^2}}

 \title{Quadratic BSDEs with jumps: a fixed-point approach\footnote{Research partly supported by the Chair {\it Financial Risks} of the {\it Risk Foundation} sponsored by Soci\'et\'e G\'en\'erale, the Chair {\it Derivatives of the Future} sponsored by the {F\'ed\'eration Bancaire Fran\c{c}aise}, and the Chair {\it Finance and Sustainable Development} sponsored by EDF and Calyon.}
}
\author{Nabil {\sc Kazi-Tani}\footnote{CMAP, Ecole Polytechnique, Paris, nabil.kazitani@polytechnique.edu.} \and Dylan {\sc Possama\"{i}}\footnote{CEREMADE, Universit\'e Paris Dauphine, possamai@ceremade.dauphine.fr.}
      \and Chao {\sc Zhou}\footnote{Department of Mathematics, National University of Singapore, Singapore, matzc@nus.edu.sg. Part of this work was carried out while the author was working at CMAP, Ecole Polytechnique,  whose financial support is kindly acknowledged.} }          
 \date{\today}

\begin{document}

 \maketitle

\vspace{3mm}

 \begin{abstract}
\vspace{10mm}

In this article, we prove the existence of bounded solutions of quadratic backward SDEs with jumps, that is to say for which the generator has quadratic growth in the variables $(z,u)$. From a technical point of view, we use a direct fixed point approach as in Tevzadze \cite{tev}, which allows us to obtain existence and uniqueness of a solution when the terminal condition is small enough. Then, thanks to a well-chosen splitting, we recover an existence result for general bounded solution. Under additional assumptions, we can obtain stability results and a comparison theorem, which as usual imply uniqueness.

\vspace{0.8em}
\noindent{\bf Key words:} BSDEs, quadratic growth, jumps, fixed-point theorem.
\vspace{5mm}

\noindent{\bf AMS 2000 subject classifications:} 60H10, 60H30
\end{abstract}
\newpage

\section{Introduction}

Motivated by duality methods and maximum principles for optimal stochastic control, Bismut studied in \cite{bis} a linear backward stochastic differential equation (BSDE). In their seminal paper \cite{pardpeng}, Pardoux and Peng generalized such equations to the non-linear Lipschitz case and proved existence and uniqueness results in a Brownian framework. Since then, a lot of attention has been given to BSDEs and their applications, not only in stochastic control, but also in theoretical economics, stochastic differential games and financial mathematics. In this context, the generalization of Backward SDEs to a setting with jumps enlarges again the scope of their applications, for instance to insurance modeling, in which jumps are inherent (see for instance Liu and Ma \cite{liu}). Li and Tang \cite{li} were the first to obtain a wellposedness result for Lipschitz BSDEs with jumps, using a fixed point approach similar to the one used in \cite{pardpeng}.

\vspace{0.3em}
Let us now precise the structure of these equations in a discontinuous setting. Given a filtered probability space $(\Omega,\mathcal F,\left\{\mathcal F_t\right\}_{0\leq t\leq T},\mathbb P)$ generated by an $\mathbb R^d$-valued Brownian motion $B$ and a random measure $\mu$ with compensator $\nu$, solving a BSDEJ with generator $g$ and terminal condition $\xi$ consists in finding a triple of progressively measurable processes $(Y,Z,U)$ such that for all $t \in [0,T]$, $\mathbb P-a.s.$
\begin{align}
 Y_t=\xi +\int_t^T g_s(Y_s,Z_s,U_s)ds-\int_t^T Z_s dB_s -\int_t^T \int_{\R^d\backslash \{0\}} U_s(x)(\mu-\nu)(ds,dx). \label{def_bsdej}
\end{align}
We refer the reader to Section \ref{notations_qbsdej} for more precise definitions and notations. In this paper, $g$ will be supposed to satisfy a Lipschitz-quadratic growth property. More precisely, $g$ will be Lipschitz in $y$, and will satisfy a quadratic growth condition in $(z,u)$ (see Assumption \ref{assump:hquad}(iii) below). The interest for such a class of quadratic BSDEs has increased a lot in the past few years, mainly due to the fact that they naturally appear in many stochastic control problems, for instance involving utility maximization (see among many others \cite{ekr} and \cite{him}).

\vspace{0.3em}
When the filtration is generated only by a Brownian motion, the existence and uniqueness of quadratic BSDEs with a bounded terminal condition has been first treated by Kobylanski \cite{kob}. Using an exponential transformation, she managed to fall back into the scope of BSDEs with a coefficient having linear growth. Then the wellposedness result for quadratic BSDEs is obtained by means of an approximation method. The main difficulty lies then in proving that the martingale part of the approximation converges in a strong sense. This result has then been extended in several directions, to a continuous
setting by Morlais \cite{morlais}, to unbounded solutions by Briand and Hu \cite{bh} or more
recently by Mocha and Westray \cite{moc}. In particular cases, several authors managed to obtain further results, to name but a few, see Hu and Schweizer \cite{hu}, Hu, Imkeller and M\"{u}ller \cite{him}, Mania and Tevzadze \cite{man} or Delbaen, et al. \cite{del}. This approach was later totally revisited by Tevzadze \cite{tev}, who gave a direct proof in the Lipschitz-quadratic setting. His methodology is fundamentally different, since he uses a fixed-point argument to obtain existence of a solution for small terminal condition, and then pastes solutions together in the general bounded case. In this regard, there is no longer any need to obtain the difficult strong convergence result needed by Kobylanski \cite{kob}. More recently, applying yet a completely different approach using now a forward point of view and stability results for a special class of quadratic semimartingales, Barrieu and El Karoui \cite{elkarbar} generalized the above results. Their approach has the merit of greatly simplifying the problem of strong convergence of the martingale part when using approximation arguments, since they rely on very general semimartingale convergence results. Notice that this approach was, partially, present in an earlier work of Cazanave, Barrieu and El Karoui \cite{elkarcaz}, but limited to a bounded framework.

\vspace{0.3em}
Nonetheless, when it comes to quadratic BSDEs in a discontinuous setting, the literature is far less abounding. Until very recently, the only existing results concerned particular cases of quadratic BSDEs, which were exactly the ones appearing in utility maximization or indifference pricing problems in a jump setting. Thus, Becherer \cite{bech} first studied bounded solutions to BSDEs with jumps in a finite activity setting, and his general results were improved by Morlais \cite{mor}, who proved existence of the solution to a special quadratic BSDE with jumps, which naturally appears in a utility maximization problem, using the same type of techniques as Kobylanski. The first breakthrough in order to tackle the general case was obtained by Ngoupeyou \cite{ngou} in his PhD thesis, and the subsequent papers by El Karoui, Matoussi and Ngoupeyou \cite{elmatn} and by Jeanblanc, Matoussi and Ngoupeyou \cite{jmn}. They non-trivially extended the techniques developed in \cite{elkarbar} to a jump setting, and managed to obtain existence of solutions for quadratic BSDEs with non-bounded terminal conditions. We emphasize that some of our arguments were inspired by their techniques and the ones developed in \cite{elkarbar}. Nonetheless, as explained throughout the paper, our approach follows a completely different direction and allows in some cases to consider BSDEs which are outside of the scope of \cite{elmatn}, even though, unlike them, we are constrained to work with bounded terminal conditions. Moreover, at least for small terminal conditions, our approach allows to obtain a wellposedness theory for multidimensional quadratic BSDEs with jumps.

\vspace{0.3em}
After the completion of this paper, we became aware of a very recent result of Laeven and Stadje \cite{laev} who proved a general existence result for BSDEJs with convex generators, using verification arguments. We emphasize that our approach is very different and do not need any convexity assumption in order to obtain existence of a solution. Nonetheless, their result and ours do not imply each other. 

\vspace{0.3em}
Our aim here is to extend the fixed-point methodology of Tevzadze \cite{tev} to the case of a discontinuous filtration. We first obtain an existence result for a terminal condition $\xi$ having a $\No{\cdot}_{\infty}$-norm which is small enough. Then the result for any $\xi$ in $\mathbb L^{\infty}$ follows by splitting $\xi$ in pieces having a small enough norm, and then pasting the obtained solutions to a single equation. Since we deal with bounded solutions, the space of BMO martingales will play a particular role in our setting. We will show that it is indeed the natural space for the continuous and the pure jump martingale terms appearing in the BSDE \ref{def_bsdej}, when $Y$ is bounded. When it comes to uniqueness of a solution in this framework with jumps, we need additional assumptions on the generator $g$ for a comparison theorem to hold. Namely, we will use on the one hand the Assumption \ref{assump.roy}, which was first introduced by Royer \cite{roy} in order to ensure the validity of a comparison theorem for Lipschitz BSDEs with jumps, and on the other hand a convexity assumption which was already considered by Briand and Hu \cite{bh2} in the continuous case. We extend here these comparison theorems to our setting (Proposition \ref{prop.comp}), and then use them to give a uniqueness result.

\vspace{0.3em}
This wellposedness result for bounded quadratic BSDEs with jumps opens the way to many possible applications. Barrieu and El Karoui \cite{elkarbar2} used quadratic BSDEs to define time consistent convex risk measures and study their properties. The extension of some of these results to the case with jumps is the object of our accompanying paper \cite{kpz4}.
 
\vspace{0.3em}
The rest of this paper is organized as follows. In Section \ref{section.1}, we give all the notations and present the natural spaces and norms in our framework. Then in Section \ref{sec.qbsdej} we provide the definition of BSDE with jump, we give the main assumptions on our generator and we prove several a priori estimates for the solution corresponding solution. Next, in Sections \ref{sec.ex1} and \ref{sec.ex2} we prove an existence result for a small enough terminal condition which we then extend to the general bounded case. Finally, Section \ref{sec.comp} is devoted to the obtention of comparison theorems and stability results for our class of BSDEJs.

\section{Preliminaries} \label{section.1}

\subsection{Notations} \label{notations_qbsdej}

We consider in all the paper a filtered probability space $\left(\Omega,\mathcal F, \left\lbrace \mathcal F_t\right\rbrace_{0\leq t\leq T},\mathbb P\right)$, whose filtration satisfies the usual hypotheses of completeness and right-continuity. We suppose that this filtration is generated by a $d$-dimensional Brownian motion $B$ and an independent integer valued random measure $\mu(\omega,dt,dx)$ defined on $\mathbb R^+\times E$, with compensator $\lambda(\omega,dt,dx)$. $\widetilde \Omega:= \Omega \times \mathbb R^+ \times E$ is equipped with the $\sigma$-field  $\widetilde{\mathcal P}:= \mathcal P \times \mathcal E$, where $\mathcal P$ denotes  the predictable $\sigma$-field on $\Omega \times \mathbb R^+$ and $\Ec$ is the Borel $\sigma$-field on $E$.

\vspace{0.3em}
To guarantee the existence of the compensator $\lambda(\omega,dt,dx)$, we assume that for each $A$ in $\Bc(E)$ and each $\omega$ in $\Omega$, the process $X_t:= \mu(\omega,A,[0,t]) \in \Ac^+_{loc}$, which means that there exists an increasing sequence of stopping times $(T_n)$ such that $T_n \to + \infty$ a.s. and the stopped processes $X_t^{T_n}$ are increasing, c\`adl\`ag, adapted and satisfy $\E[X_{\infty}]<+\infty$.

\vspace{0.3em}
We assume in all the paper that $\lambda$ is absolutely continuous with respect to the Lebesgue measure $dt$, i.e. $\lambda(\omega,dt,dx)=\nu_t(\omega,dx)dt$. Finally, we denote $\widetilde\mu$ the compensated jump measure
$$\widetilde\mu(\omega,dx,dt) = \mu(\omega,dx,dt) - \nu_t(\omega,dx)\, dt.$$


In our setting, we emphasise that  we allow the compensator of the jump measure to be a random measure unlike most of the literature where it is a classic L\'evy measure (see however \cite{bech} for a similar approach). This will not increase the complexity of our proofs, provided that the martingale representation property of Assumption \ref{martingale_representation} below holds true.

\vspace{0.3em}
Following Tang and Li \cite{li} and Barles et al. \cite{barles}, the definition of a BSDE with jumps is then 
\begin{Definition}\label{def_bsdej2}
Let $\xi$ be a $\mathcal F_T$-measurable random variable. A solution to the BSDEJ with terminal condition $\xi$ and generator $g$ is a triple $(Y,Z,U)$ of progressively measurable processes such that
\begin{equation}
Y_t=\xi+\int_t^Tg_s(Y_s,Z_s,U_s)ds-\int_t^TZ_sdB_s-\int_t^T\int_{E} U_s(x)\widetilde\mu(dx,ds),\ t\in[0,T],\ \mathbb P-a.s.
\label{eq:bsdej}
\end{equation}
\end{Definition}
where $g:\Omega\times[0,T]\times\mathbb R\times\mathbb R^d\times \Ac(E) \rightarrow \mathbb R$ is a given application and 
$$\Ac(E):=\left\{u: \, E \rightarrow \R,\, \Bc(E)-\text{measurable} \right\}.$$

Then, the processes $Z$ and $U$ are supposed to satisfy the minimal assumptions so that the quantities in \reff{eq:bsdej} are well defined, namely $(Z,U) \in \mathcal Z \x \mathcal U$, where $\mathcal Z$ (resp. $\Uc$) denotes the space of all $\mathbb F$-predictable $\mathbb R^d$-valued processes $Z$ (resp. $\mathbb F$-predictable functions $U$) with
$$\int_0^T\abs{Z_t}^2dt <+\infty,\ \left(\text{resp. }\int_0^T\int_{E}\abs{U_t(x)}^2 \nu_t(dx)dt <+\infty\right),\ \mathbb P-a.s.$$

\begin{Remark}
Notice that this is a particular case of the framework considered earlier by El Karoui and Huang \cite{elkh} and El Karoui et al. \cite{elkaroui} where the filtration is more general, and therefore they do not have an explicit form for the martingale part which is orthogonal to the Brownian one. Here, knowing explicitly this orthogonal martingale allows to have a dependance of the generator in it, through the predictable function $U$. Here $U$ plays a role analogous to the quadratic variation in the continuous case. However, there are some notable differences, since for each $t$, $U_t$ is a function mapping $E$ to $\R$. This is why the treatment of the dependence in $u$ in the assumptions for the generator is not symmetric to the treatment of the dependence in $z$, and in particular we deal with Fr\'echet derivatives with respect to $u$ (see Assumption \ref{lipschitz_assumption} for more precise statements).
\end{Remark}

\subsection{Standard spaces and norms}

We introduce the following norms and spaces for any $p\geq 1$.

\vspace{0.3em}

$\mathcal S^\infty$ is the space of $\mathbb R$-valued c\`adl\`ag and $\mathcal F_t$-progressively measurable processes $Y$ such that
$$\No{Y}_{\mathcal S^\infty}:=\underset{0\leq t\leq T}{\sup}\No{Y_t}_\infty<+\infty.$$

$\mathbb H^p$ is the space of $\mathbb R^d$-valued and $\mathcal F_t$-progressively measurable processes $Z$ such that
$$\No{Z}^p_{\mathbb H^p}:=\mathbb E\left[\left(\int_0^T\abs{Z_t}^2dt\right)^{\frac p2}\right]<+\infty.$$

The two spaces above are the classical ones in the BSDE theory in continuous filtrations. We introduce finally a space which is specific to the jump case, and which plays the same role for $U$ as $\mathbb H^p$ for $Z$. $\mathbb J^p$ is the space of predictable and $\mathcal E$-measurable applications $U:\Omega\times[0,T]\times E$ such that
$$\No{U}^p_{\mathbb J^p}:=\mathbb E\left[\left(\int_0^T\int_E\abs{U_s(x)}^2\nu_s(dx)ds\right)^{\frac p2}\right]<+\infty.$$

\subsection{A word on c\`adl\`ag BMO martingales}

The recent literature on quadratic BSDEs is very rich on remarks and comments about the deep theory of continuous BMO martingales. However, it is clearly not as well documented when it comes to c\`adl\`ag BMO martingales, whose properties are crucial in this paper. Indeed, apart from some remarks in the book by Kazamaki \cite{kaz}, the extension to the c\`adl\`ag case of the classical results of BMO theory, cannot always be easily found. Our main goal in this short subsection is to give a rapid overview of the existing literature and results concerning BMO martingales with c\`adl\`ag trajectories, with an emphasis where the results differ from the continuous case. Let us start by recalling some notations and definitions.

\vspace{0.3em}
$\rm{BMO}$ is the space of square integrable c\`adl\`ag $\mathbb R^d$-valued martingales $M$ such that
$$\No{M}_{\rm{BMO}}:=\underset{\tau\in\mathcal T_0^T}{\esup^\mathbb P}\No{\mathbb E_\tau\left[\left(M_T-M_{\tau^-}\right)^2\right]}_{\infty}<+\infty,$$
where for any $t\in[0,T]$, $\mathcal T_t^T$ is the set of $(\mathcal F_s)_{0\leq s\leq T}$-stopping times taking their values in $[t,T]$.

\vspace{0.3em}
$\mathbb J^2_{\rm{BMO}}$ is the space of predictable and $\mathcal E$-measurable applications $U:\Omega\times[0,T]\times E$ such that
$$\No{U}^2_{\mathbb J^2_{\rm{BMO}}}:=\No{\int_0^.\int_EU_s(x)\widetilde\mu(dx,ds)}_{ \rm{BMO}}<+\infty.$$

\vspace{0.3em}
$\mathbb H^2_{\rm{BMO}}$ is the space of $\mathbb R^d$-valued and $\mathcal F_t$-progressively measurable processes $Z$ such that
$$\No{Z}^2_{\mathbb H^2_{\rm{BMO}}}:=\No{\int_0^.Z_sdB_s}_{\rm{BMO}}<+\infty.$$

\vspace{0.3em}
As soon as the process $\langle M\rangle$ is defined for a martingale $M$, which is the case if for instance $M$ is locally square integrable, then it is easy to see that $M\in  \rm{BMO}$ if the jumps of $M$ are uniformly bounded in $t$ by some positive constant $C$ and
$$\underset{\tau\in\mathcal T_0^T}{\esup^\mathbb P}\No{\mathbb E_\tau\left[\langle M\rangle_T-\langle M\rangle_{\tau}\right]}_{\infty}\leq C.$$
Furthermore the BMO norm of $M$ is then smaller than $2C$. We also recall the so called energy inequalities (see \cite{kaz} and the references therein). Let $Z\in\mathbb H^2_{\rm{BMO}}$, $U\in\mathbb J^2_{\rm{BMO}}$ and $p\geq 1$. Then we have
\begin{align}\label{energy}
\nonumber&\mathbb E\left[\left(\int_0^T\abs{Z_s}^2ds\right)^p\right]\leq 2p!\left(4\No{Z}_{\mathbb H^2_{\rm{BMO}}}^2\right)^p\\
&\mathbb E\left[\left(\int_0^T\int_EU_s^2(x)\nu_s(dx)ds\right)^p\right]\leq 2p!\left(4\No{U}_{\mathbb J^2_{\rm{BMO}}}^2\right)^p.
\end{align}

Let us now turn to more precise properties and estimates for BMO martingales. It is a classical result (see \cite{kaz}) that the Dol\'eans-Dade exponential of a continuous BMO martingale is a uniformly integrable martingale. Things become a bit more complicated in the c\`adl\`ag case, and more assumptions are needed. Let us first define the Dol\'eans-Dade exponential of a square integrable martingale $X$, denoted $\mathcal E(X)$. This is as usual the unique solution $Z$ of the SDE
$$Z_t=1+\int_0^tZ_{s^-}dX_s, \ \mathbb P-a.s.,$$
and is given by the formula
$$\mathcal E(X)_t=e^{X_t-\frac12<X^c>_t}\prod_{0<s\leq t}(1+\Delta X_s)e^{-\Delta X_s},\ \mathbb P-a.s.$$

One of the first results concerning Dol\'eans-Dade exponential of BMO martingales was proved by Dol\'eans-Dade and Meyer \cite{dolm}. They showed that

\begin{Proposition}
Let $M$ be a c\`adl\`ag BMO martingale such that $\No{M}_{\rm{BMO}}<1/8$. Then $\mathcal E(M)$ is a strictly positive uniformly integrable martingale.
\end{Proposition}

The constraint on the norm of the martingale being rather limiting for applications, this result was subsequently improved by Kazamaki \cite{kaz2}, where the constraints is now on the jumps of the martingale  

\begin{Proposition}
Let $M$ be a c\`adl\`ag BMO martingale such that there exists $\delta >0$ with $\Delta M_t\geq -1+\delta$, for all $t\in [0,T]$, $\mathbb P-a.s.$ Then $\mathcal E(M)$ is a strictly positive uniformly integrable martingale.
\end{Proposition}

Furthermore, we emphasize, as recalled in the counter-example of Remark $2.3$ in \cite{kaz}, that a complete generalization to the c\`adl\`ag case is not possible. We also refer the reader to L\'epingle and M\'emin \cite{lepinmemin1} and \cite{lepinmemin2} for general sufficient conditions for the uniform integrability of Dol\'eans-Dade exponentials of c\`adl\`ag martingales. This also allows us to obtain immediately a Girsanov Theorem in this setting, which will be extremely useful throughout the paper.

\begin{Proposition}\label{girsanov}
Let us consider the following c\`adl\`ag martingale $M$ 
$$M_t:=\int_0^t\varphi_sdB_s+\int_0^t\int_E\gamma_s(x)\widetilde\mu(dx,ds),\ \mathbb P-a.s.,$$
where $\gamma$ is bounded and $(\varphi,\gamma)\in\mathbb H^2_{\rm{BMO}}\times\mathbb J^2_{\rm{BMO}}$ and where there exists $\delta >0$ with $\gamma_t\geq -1+\delta$, $\mathbb P\times d\nu_t-a.e.$, for all $t\in [0,T]$.

\vspace{0.3em}
Then, the probability measure $\mathbb Q$ defined by $\frac{d\mathbb Q}{d\mathbb P}=\mathcal E\left(M_.\right),$ is indeed well-defined and starting from any $\mathbb P$-martingale, by, as usual, changing adequately the drift and the jump intensity, we can obtain a $\mathbb Q$-martingale.
\end{Proposition}

We now address the question of the so-called reverse H\"older inequality, which implies in the continuous case that if $M$ is a BMO martingale, there exists some $r>1$ such that $\mathcal E(M)$ is $L^r$-integrable. As for the previous result on uniform integrability, this was extended to the c\`adl\`ag case first in \cite{dolm} and \cite{kaz3}, with the additional assumption that the BMO norm or the jumps of $M$ are sufficiently small. The following generalization is taken from \cite{izu}
\begin{Proposition}\label{prop.expor}
Let $M$ be a c\`adl\`ag BMO martingale such that there exists $\delta >0$ with $\Delta M_t\geq -1+\delta$, for all $t\in [0,T]$, $\mathbb P-a.s.$ Then $\mathcal E(M)$ is in $L^r$ for some $r>1$.
\end{Proposition}

%

\section{Quadratic BSDEs with jumps}\label{sec.qbsdej}
\subsection{The non-linear generator}

Following the Definition \ref{def_bsdej} of BSDEs with jumps, we need now to specify in more details the assumptions we make on the generator $g$. The most important one in our setting will be the quadratic growth assumption of Assumption \ref{assump:hquad}(ii) below. It is the natural generalization to the jump case of the usual quadratic growth assumption in $z$. Before proceeding further, let us define the following function
$$j_t(u):=\int_E\left(e^{u(x)}-1-u(x)\right)\nu_t(dx).$$

This function $j(u)$ plays the same role for the variable $u$ as the square function for the variable $z$. In order to understand this, let us consider the following "simplest" quadratic BSDE with jumps
$$y_t=\xi+\int_t^T\left(\frac\gamma2\abs{z_s}^2+\frac1\gamma j_s(\gamma u_s)\right)ds-\int_t^Tz_sdB_s-\int_t^T\int_Eu_s(x)\widetilde\mu(dx,ds), \ t\in[0,T],\\ \mathbb P-a.s.$$

Then a simple application of It\^o's formula gives formally
$$e^{\gamma y_t}=e^{\gamma \xi}-\int_t^Te^{\gamma y_s}z_sdB_s-\int_t^T\int_Ee^{\gamma y_{s^-}}\left(e^{\gamma u_s(x)}-1\right)\widetilde\mu(dx,ds),\ t\in[0,T],\ \mathbb P-a.s.$$

Still formally, taking the conditional expectation above gives finally
$$y_t=\frac1\gamma\ln\left(\mathbb E_t\left[e^{\gamma \xi}\right]\right),\ t\in[0,T],\ \mathbb P-a.s.,$$
and we recover the so-called entropic risk measure which in the continuous case corresponds to a BSDE with generator $\frac{\gamma}{2}\abs{z}^2$.

\vspace{0.8em}
Of course, for the above to make sense, the function $j$ must at the very least be well defined. A simple application of Taylor's inequalities shows that if the function $x\mapsto u(x)$ is bounded $d\nu_t-a.e.$ for every $0\leq t\leq T$, then we have for some constant $C>0$
$$0\leq e^{u(x)}-1-u(x)\leq Cu^2(x), \text{ $d\nu_t-a.e.$ for every $0\leq t\leq T$.}$$

Hence, if we introduce for $1<p\leq +\infty$ the spaces
$$L^p(\nu):=\left\{u,\ \text{$\mathcal E$-measurable, such that $u\in L^p(\nu_t)$ for all $0\leq t\leq T$}\right\},$$
then $j$ is well defined on $L^2(\nu)\cap L^\infty(\nu)$. We now give our quadratic growth assumption on $g$.

\begin{Assumption}\label{assump:hquad}[Quadratic growth]

\vspace{0.3em}
\rm{(i)} For fixed $(y,z,u)$, $g$ is $\mathbb{F}$-progressively measurable.

\vspace{0.3em}
\rm{(ii)} For any $p\geq 1$
\begin{equation}\label{inte}
\underset{\tau\in\mathcal T_0^T}{\esup^\mathbb P}\ \mathbb E_\tau\left[\left(\int_\tau^T\abs{g_t(0,0,0)}dt\right)^p\right]<+\infty, \ \mathbb P-a.s.
\end{equation} 

\rm{(iii)} $g$ has the following growth property. There exists $(\beta,\gamma)\in \mathbb R_+\times \mathbb R^*_+$ and a positive predictable process $\alpha$ satisfying the same integrability condition \reff{inte} as $g_t(0,0,0)$, such that for all $(\omega,t,y,z,u)$
\begin{align}
 -\alpha_t-\beta\abs{y}-\frac\gamma2\abs{z}^2-\frac{j_t(-\gamma u)}{\gamma} \leq g_t(\omega,y,z,u)-g_t(0,0,0)\leq \alpha_t+\beta\abs{y}+\frac\gamma2\abs{z}^2+\frac{j_t(\gamma u)}{\gamma} .\label{eq_quadratique}
\end{align}
\end{Assumption}

\begin{Remark}\label{remrem}
We emphasize that unlike the usual quadratic growth assumptions for continuous BSDEs, condition \reff{eq_quadratique} is not symmetric. It is mainly due to the fact that unlike the functions $\abs{.}$ and $\abs{.}^2$, the function $j$ is not even. Moreover, with this non-symmetric condition, it is easily seen that if $Y$ is a solution to equation \reff{eq:bsdej} with a generator satisfying the condition \reff{eq_quadratique}, then $-Y$ is also a solution to a BSDE whose generator satisfy the same condition \reff{eq_quadratique}. More precisely, if $(Y,Z,U)$ solves equation \reff{eq:bsdej}, then $(-Y,-Z,-U)$ solves the BSDEJ with terminal condition $-\xi$ and generator $\widetilde g_t(y,z,u):=-g_t(-y,-z,-u)$ which clearly also satisfies \reff{eq_quadratique}. This will be important for the proof of Lemma \ref{lemma.bmo}.

\vspace{0.3em}
We also want to insist on the structure which appears in \reff{eq_quadratique}. Indeed, the constant $\gamma$ in front of the quadratic term in $z$ is the same as the one appearing in the term involving the function $j$. As already seen for the entropic risk measure above, if the constants had been different, say respectively $\gamma_1$ and $\gamma_2$,  the exponential transformation would have failed. Moreover, since the function $\gamma\mapsto \gamma^{-1}j_t(\gamma u)$ is not monotone, then we cannot increase or decrease $\gamma_1$ and $\gamma_2$ to recover the desired estimate \reff{eq_quadratique}. Moreover, we emphasize that such a structure already appeared in \cite{elmatn}, \cite{jmn} and \cite{ngou}, where it was also crucial in order to obtain existence. Notice however that thanks to our particular context of bounded terminal conditions, we will show that in some cases, we are no longer constrained by this structure (see Remark \ref{rem.assumptions}).
\end{Remark}

\subsection{First a priori estimates for the solution}
We first prove a result showing a link between the BMO spaces and quadratic BSDEs with jumps, a property which is very well known in the continuous case since the paper by Hu, Imkeller and M\"{u}ller \cite{him}, and which also appears in \cite{mor} and \cite{ngou}. We emphasize that only Assumption \ref{assump:hquad} is necessary to obtain it. Before proceeding, we define for every $x\in\mathbb R$ and every $\eta\neq0$, $h_\eta(x):=(e^{\eta x}-1-\eta x)/\eta$. The function $h_\eta$ already appears in our growth Assumption \ref{assump:hquad}(ii), and the following trivial property that it satisfies is going to be crucial for us
\begin{equation}
h_{2\eta}(x)=\frac{\left(e^{\eta x}-1\right)^2}{2\eta}+ h_\eta(x).
\label{eq:oulala}
\end{equation}

We also give the two following inequalities which are of the utmost importance in our jump setting. We emphasize that the first one is trivial, while the second one can be proved using simple but tedious algebra.
\begin{align}\label{1}
2\leq e^x+e^{-x},\text{ $\forall x\in\mathbb R$},\ \ x^2\leq a\left(e^x-1\right)^2+ \frac{\left(1-e^{-x}\right)^2}{a},\text{ $\forall(a,x)\in\mathbb R^*_+\times\mathbb R$.}
\end{align}

We then have the following Lemma (which is closely related to Proposition $8$ in \cite{ngou}). 
\begin{Lemma}\label{lemma.bmo}
Let Assumption \ref{assump:hquad} hold. Assume that $(Y,Z,U)$ is a solution to the BSDEJ \reff{eq:bsdej} such that $(Z,U)\in \mathcal Z\times \mathcal U$, the jumps of $Y$ are bounded and
\begin{equation}
\underset{\tau\in\mathcal T_0^T}{\esup^{\P}}\ \mathbb E_\tau\left[\exp\left(2\gamma\underset{\tau\leq t\leq T}{\sup}\pm Y_t\right)\vee\exp\left(4\gamma\underset{\tau\leq t\leq T}{\sup}\pm Y_t\right)\right]<+\infty, \; \mathbb P-a.s.
\label{eq:3}
\end{equation}
 
\vspace{0.3em}
Then $Z\in\mathbb H^2_{\rm{BMO}}$ and $U\in \mathbb J^2_{\rm{BMO}}\cap L^\infty(\nu)$.
\end{Lemma}

\vspace{0.3em}
\proof
First of all, since the size of the jumps of $Y$ is bounded, there exists a version of $U$, that is to say that there exists a predictable function $\widetilde U$ such that for all $t\in[0,T]$
$$\int_E\abs{\widetilde U_t(x)-U_t(x)}^2\nu_t(dx)=0,\ \mathbb P-a.s.,$$
and such that $|\widetilde{U}_t(x)|\leq C,\text{ for all $t$, $\mathbb P-a.s$}.$ For the sake of simplicity, we will always consider this version and we still denote it $U$. For the proof of this result, we refer to Morlais \cite{mor}.

\vspace{0.3em}
Let us consider the following processes
$$\int_0^Te^{2\gamma Y_t}Z_tdB_t\text{ and }\int_0^Te^{2\gamma Y_{t^-}}\left(e^{2\gamma U_t(x)}-1\right)\widetilde\mu(dx,dt).$$

We will first show that they are local martingales. Indeed, we have
\begin{align*}
\int_0^Te^{4\gamma Y_t}Z_t^2dt\leq \exp\left(4\gamma \underset{0\leq t\leq T}{\sup} Y_t\right)\int_0^TZ_t^2dt<+\infty, \ \mathbb P-a.s., 
\end{align*}
since $Z\in\mathcal Z$ and \reff{eq:3} holds. Similarly, we have
\begin{align*}
\int_0^T\int_Ee^{4\gamma Y_t}U_t^2(x)\nu_t(dx)dt\leq \exp\left(4\gamma \underset{0\leq t\leq T}{\sup} Y_t\right)\int_0^T\int_EU_t^2(x)\nu_t(dx)dt<+\infty, \ \mathbb P-a.s., 
\end{align*}
since $U\in\mathcal U$ and \reff{eq:3} holds.

\vspace{0.3em}
Let now $(\tau_n)_{n\geq 1}$ be a localizing sequence for the $\mathbb P$-local martingales above. By It\^o's formula under $\mathbb P$ applied to $e^{2\gamma Y_t}$, we have for every $\tau \in \mathcal T^T_0$
\begin{align*}
&\frac{4\gamma^2}{2}\int_\tau^{\tau_n}e^{2\gamma Y_t}\abs{Z_t}^2dt+2\gamma\int_\tau^{\tau_n}\int_Ee^{2\gamma Y_t}h_{2\gamma}\left(U_t(x)\right)\nu_t(dx)dt\\
&=e^{2\gamma Y_{\tau_n}}-e^{2\gamma Y_{\tau}}+2\gamma\int_\tau^{\tau_n}e^{2\gamma Y_t}g_t(Y_t,Z_t,U_t)dt-2\gamma\int_\tau^{\tau_n}e^{2\gamma Y_{t}}Z_tdB_t\\
&\hspace{0.9em}-2\gamma\int_\tau^{\tau_n}\int_Ee^{2\gamma Y_{t^{-}}}\left(e^{2\gamma U_t(x)}-1\right)\widetilde\mu(dx,dt)\\
&\leq e^{2\gamma Y_{\tau_n}}-e^{2\gamma Y_{\tau}}+2\gamma\int_\tau^{\tau_n}e^{2\gamma Y_t}\left(\alpha_t+\abs{g_t(0,0,0)}+\beta\abs{Y_t}\right)dt\\
&\hspace{0.9em}+2\gamma\int_\tau^{\tau_n}e^{2\gamma Y_t}\left(\frac{\gamma}{2}\abs{Z_t}^2+\int_Eh_\gamma\left(U_t(x)\right)\nu_t(dx)\right)dt-2\gamma\int_\tau^{\tau_n}e^{2\gamma Y_{t}}Z_tdB_t\\
&\hspace{0.9em}-2\gamma\int_\tau^{\tau_n}\int_Ee^{2\gamma Y_{t^{-}}}\left(e^{2\gamma U_t(x)}-1\right)\widetilde\mu(dx,dt),\ \mathbb P-a.s.
\end{align*}

\vspace{0.3em}
Now the situation is going to be different from the continuous case, and the property \reff{eq:oulala} is going to be important. Indeed, we can take conditional expectation and thus obtain
\begin{align*}
&\mathbb E_\tau\left[\gamma^2\int_\tau^{\tau_n}e^{2\gamma Y_t}\abs{Z_t}^2dt+\int_\tau^{\tau_n}\int_Ee^{2\gamma Y_t}\left(e^{\gamma U_t(x)}-1\right)^2\nu_t(dx)dt\right]\\
&\leq C\left(1+\mathbb E_\tau\left[\left(\int_\tau^{\tau_n}\left(\alpha_t+\abs{g_t(0,0,0)}\right)dt\right)^2+\exp\left(2\gamma \underset{\tau\leq t\leq T}{\sup}Y_{t}\right) + \exp\left(4\gamma \underset{\tau\leq t\leq T}{\sup}Y_{t}\right)\right]\right)\\
&\leq C\left(1+\mathbb E_\tau\left[\exp\left(2\gamma \underset{\tau\leq t\leq T}{\sup}Y_{t}\right)\vee \exp\left(4\gamma \underset{\tau\leq t\leq T}{\sup}Y_{t}\right)\right]\right),
\end{align*}
where we used the inequality $2ab\leq a^2+b^2$, the fact that for all $x\in\mathbb R$, $\abs{x}e^x\leq C(1+e^{2x})$ for some constant $C>0$ (which as usual can change value from line to line) and the fact that Assumption \ref{assump:hquad}(ii) and (iii) hold.

\vspace{0.3em}
Using Fatou's lemma and the monotone convergence theorem, we obtain
\begin{align}\label{eq:bmo1}
\nonumber&\mathbb E_\tau\left[\gamma^2\int_\tau^{T}e^{2\gamma Y_t}\abs{Z_t}^2dt+\int_\tau^{T}\int_Ee^{2\gamma Y_t}\left(e^{\gamma U_t(x)}-1\right)^2\nu_t(dx)dt\right]\\
&\leq C\left(1+\underset{\tau\in\mathcal T_0^T}{\esup^\mathbb P}\ \mathbb E_\tau\left[ \exp\left(2\gamma \underset{\tau\leq t\leq T}{\sup}Y_{t}\right)\vee \exp\left(4\gamma \underset{\tau\leq t\leq T}{\sup}Y_{t}\right)\right]\right).
\end{align}

Now, we apply the above estimate for the solution $(-Y,-Z,-U)$ of the BSDEJ with terminal condition $-\xi$ and generator $\widetilde g_t(y,z,u):=-g_t(-y,-z,-u)$, which still satisfies Assumption \ref{assump:hquad} (see Remark \ref{remrem})
\begin{align}\label{eq:bmo2}
\nonumber&\mathbb E_\tau\left[\gamma^2\int_\tau^{T}e^{-2\gamma Y_t}\abs{Z_t}^2dt+\int_\tau^{T}\int_Ee^{-2\gamma Y_t}\left(e^{-\gamma U_t(x)}-1\right)^2\nu_t(dx)dt\right]\\
&\leq C\left(1+\underset{\tau\in\mathcal T_0^T}{\esup^\mathbb P}\ \mathbb E_\tau\left[ \exp\left(2\gamma \underset{\tau\leq t\leq T}{\sup}\left(-Y_{t}\right)\right)\vee \exp\left(4\gamma \underset{\tau\leq t\leq T}{\sup}\left(-Y_{t}\right)\right)\right]\right).
\end{align}

Let us now sum the inequalities \reff{eq:bmo1} and \reff{eq:bmo2}. We obtain
\begin{align*}
\nonumber&\mathbb E_\tau\left[\int_\tau^{T}\left(e^{2\gamma Y_t}+e^{-2\gamma Y_t}\right)\abs{Z_t}^2+\int_Ee^{2\gamma Y_t}\left(e^{\gamma U_t(x)}-1\right)^2+e^{-2\gamma Y_t}\left(e^{-\gamma U_t(x)}-1\right)^2\nu_t(dx)dt\right]\\
&\leq C\left(1+\underset{\tau\in\mathcal T_0^T}{\esup^\mathbb P}\ \mathbb E_\tau\left[ \underset{\tau\leq t\leq T}{\sup}\left\{e^{2\gamma Y_{t}}\vee e^{4\gamma Y_{t}}+e^{2\gamma \left(-Y_{t}\right)}\vee e^{4\gamma \left(-Y_{t}\right)}\right\}\right]\right).
\end{align*}

Finally, from the inequalities in \reff{1}, this shows the desired result.
\ep

\begin{Remark}
In the above Proposition, if we only assume that 
$$\mathbb E\left[\exp\left(2\gamma\underset{0\leq t\leq T}{\sup}\pm Y_t\right)\vee \exp\left(4\gamma\underset{0\leq t\leq T}{\sup}\pm Y_t\right)\right]<+\infty,$$
then the exact same proof would show that $(Z,U)\in\mathbb H^2\times\mathbb J^2$. Moreover, using the Neveu-Garsia Lemma in the same spirit as \cite{elkarbar}, we could also show that $(Z,U)\in\mathbb H^p\times\mathbb J^p$ for all $p>1$.
\end{Remark}

\vspace{0.3em}
We emphasize that the results of this Proposition highlight the fact that we do not necessarily need to consider solutions with a bounded $Y$ in the quadratic case to obtain {\it a priori} estimates. It is enough to assume the existence of some exponential moments. This is exactly the framework developed in \cite{bh} and \cite{elkarbar} in the continuous case and in \cite{elmatn} and \cite{ngou} in the jump case. It implies furthermore that it s not necessary to let the BMO spaces play a particular role in the general theory. Nonetheless, our proof of existence will rely heavily on BMO properties of the solution, and the simplest condition to obtain the estimate \reff{eq:3} is to assume that $Y$ is indeed bounded. The aim of the following Proposition is to show that we can control the $\mathcal S^\infty$ norm of $Y$ by the $L^\infty$ norm of $\xi$. Since the proof is very similar to the proof of Lemma $1$ in \cite{bh}, we will omit it.

\begin{Proposition}\label{prop.estim}
Let $\xi\in\mathbb L^\infty$. Let Assumption \ref{assump:hquad} hold and assume that $$\abs{g(0,0,0)}+\alpha\leq M,$$ for some constant $M>0$. Let $(Y,Z,U)\in \mathcal S^{\infty}\times\mathbb H^{2}\times \mathbb J^2$ be a solution to the BSDEJ \reff{eq:bsdej}. Then we have
\begin{equation*}
\abs{Y_t}\leq \gamma M\frac{e^{\beta(T-t)}-1}{\beta}+\gamma e^{\beta(T-t)}\No{\xi}_{\mathbb L^\infty},\ \mathbb P-a.s.
\end{equation*}
\end{Proposition}

\section{Existence and uniqueness for a small terminal condition}\label{sec.ex1}

The aim of this Section is to obtain an existence and uniqueness result for BSDEJs with quadratic growth when the terminal condition is small enough. However, we will need more assumptions for our proof to work. First, we assume from now on that we have the following martingale representation property. We need this assumption since we will rely on the existence results in \cite{barles} or \cite{li} which need the martingale representation.

\begin{Assumption}\label{martingale_representation}
Any local martingale $M$ with respect to the filtration $(\mathcal F_t)_{0\leq t\leq T}$ has the predictable representation property, that is to say that there exist a unique predictable process $H$ and a unique predictable function $U$ such that $(H,U)\in\mathcal Z\times\mathcal U$ and
$$M_t=M_0+\int_0^tH_sdB_s+\int_0^t\int_EU_s(x)\widetilde\mu(dx,ds), \; \mathbb P-a.s.$$
\end{Assumption}

\begin{Remark}
This martingale representation property holds for instance when the compensator $\nu$ does not depend on $\omega$, i.e when $\nu$ is the compensator of the counting measure of an additive process in the sense of Sato \cite{sato}. It also holds when $\nu$ has the particular form described in \cite{kpz1}, in which case $\nu$ depends on $\omega$.
\end{Remark} 

Of course, we also need to assume more properties for our generator $g$. Before stating them, let us describe the underlying intuitions. We want to obtain existence through a fixed point argument, therefore we have to assume some kind of control in $(y,z,u)$ of our generator. In the classical setting of \cite{pardpeng} and \cite{elkaroui}, the required contraction is obtained by using the Lipschitz property of the generator $g$ and by considering well-chosen weighted norms. More precisely, and abusing notations, they consider for some constant $\upsilon$ the spaces $\mathbb H^2_\upsilon$ consisting of progressively measurable processes $X$ such that 
$$\No{X}^2_{\mathbb H^2_\upsilon}:=\mathbb E\left[\int_0^Te^{\upsilon s}\abs{X_s}^2ds\right]<+\infty.$$

Then by choosing $\upsilon$ large enough, they can obtain a contraction in these spaces. In our context and in the context of \cite{tev}, the Lipschitz assumption for the generator, which would imply linear growth, is replaced by some kind of local Lipschitz assumption with quadratic growth. In return, it becomes generally impossible to recover a contraction. Indeed, as we will see later on, the application for which we want to find a fixed point is no longer Lipschitz but only locally Lipschitz. In these regard, it is useless for us to use weighted norms, since they can only diminish the constants intervening in our estimates. The idea is then to localize the procedure in a ball, so that the application will become Lipschitz, and then to choose the radius of this ball sufficiently small so that we actually recover a contraction. The crucial contribution of Tevzadze \cite{tev} to this problem is to show that such controls can be obtained by taking a terminal condition small enough.

\vspace{0.35em}
We now state our assumptions and refer the reader to Remark \ref{rem.assumptions} for more discussions.
\begin{Assumption}\label{lipschitz_assumption}[Lipschitz assumption]

\vspace{0.3em}
Let Assumption \ref{assump:hquad}(i),(ii) hold and assume furthermore that

\vspace{0.3em}
\rm{(i)} $g$ is uniformly Lipschitz in $y$.
$$\abs{g_t(\omega,y,z,u)-g_t(\omega,y',z,u)}\leq C\abs{y-y'}\text{ for all }(\omega,t,y,y',z,u).$$

\rm{(ii)} $\exists$ $\mu>0$ and $\phi\in\mathbb H^2_{\rm{BMO}}$ such that for all $(t,y,z,z',u)$
$$\abs{ g_t(\omega,y,z,u)- g_t(\omega,y,z',u)-\phi_t.(z-z')}\leq \mu \abs{z-z'}\left(\abs{z}+\abs{z'}\right).$$

 \rm{(iii)} $\exists$ $\mu>0$ and $\psi\in\mathbb J^2_{\rm{BMO}}$ such that for all $(t,x)$
$$C_1(1\wedge\abs{x})\leq\psi_t(x)\leq C_2(1\wedge\abs{x}),$$ where $C_2>0$, $C_1\geq-1+\delta$ where $\delta>0$. Moreover, for all $(\omega,t,y,z,u,u')$
$$\abs{ g_t(\omega,y,z,u)- g_t(\omega,y,z,u')-\langle\psi_t,u-u'\rangle_t}\leq \mu \No{u-u'}_{L^2(\nu_t)}\left(\No{u}_{L^2(\nu_t)}+\No{u'}_{L^2(\nu_t)}\right),$$
where $\langle u_1,u_2\rangle_t:=\int_Eu_1(x)u_2(x)\nu_t(dx)$ is the scalar product in $L^2(\nu_t)$.
\end{Assumption}

\vspace{0.3em}
\begin{Remark}\label{rem.assumptions}
Let us comment on the above assumptions. The first one concerning Lipschitz continuity in the variable $y$ is classical in the BSDE theory. The two others may seem a bit complicated, but as already mentioned above, they are almost equivalent to saying that the function $g$ is locally Lipschitz in $z$ and $u$. In the case of the variable $z$ for instance, those two properties would be equivalent if the process $\phi$ were bounded. Here we allow something a bit more general by letting $\phi$ be unbounded but in $\mathbb H^2_{\rm{BMO}}$. Once again, since these assumptions allow us to apply the Girsanov property of Proposition \ref{girsanov}, we do not need to bound the processes and BMO type conditions are sufficient. Moreover, Assumption \ref{lipschitz_assumption} also implies a weaker version of Assumption \ref{assump:hquad}. Indeed, it implies clearly that
$$\abs{g_t(y,z,u)-g_t(0,0,0)-\phi_t.z-\langle \psi_t,u\rangle_t}\leq C\abs{y}+\mu\left(\abs{z}^2+\No{u}^2_{L^2(\nu_t)}\right).$$

Then, for any $u\in L^2(\nu)\cap L^\infty(\nu)$ and for any $\gamma>0$, we have using the mean value Theorem
$$\frac\gamma2 e^{-\gamma\No{u}_{L^\infty(\nu)}}\No{u}^2_{L^2(\nu_t)}\leq \frac1\gamma j_t(\pm\gamma u)\leq \frac\gamma2 e^{\gamma\No{u}_{L^\infty(\nu)}}\No{u}^2_{L^2(\nu_t)}.$$

Denote $\delta g_t:=g_t(y,z,u)-g_t(0,0,0)$. We deduce using the Cauchy-Schwarz inequality and the trivial inequality $2ab\leq a^2+b^2$
\begin{align*}
&\delta g_t\leq\frac{\abs{\phi_t}^2}{2}+\frac{\No{\psi_t}^2_{L^2(\nu_t)}}{2}+C\abs{y}+\left(\mu+\frac12\right)\left(\abs{z}^2+\frac{2e^{\gamma\No{u}_{L^\infty(\nu)}}}{\gamma^2}j_t(\gamma u)\right)\\
&\delta g_t\geq-\frac{\abs{\phi_t}^2}{2}-\frac{\No{\psi_t}^2_{L^2(\nu_t)}}{2}-C\abs{y}-\left(\mu+\frac12\right)\left(\abs{z}^2+\frac{2e^{\gamma\No{u}_{L^\infty(\nu)}}}{\gamma^2}j_t(-\gamma u)\right).
\end{align*}

It is easy to check, using the energy inequalities \reff{energy} and the definition of the essential supremum,  that the term $\abs{\phi_t}^2+\No{\psi_t}^2_{L^2(\nu_t)}$ above satisfies the integrability condition \reff{inte}. Hence, we have obtained a growth property which is similar to \reff{eq_quadratique}, the only difference being that the constants appearing in the quadratic term in $z$ and the term involving the function $j$ are not the same. We thus are no longer constrained by the structure already mentioned in Remark \ref{remrem}.
\end{Remark}

We now show that if we can solve the BSDEJ \reff{eq:bsdej} for a generator $g$ satisfying Assumption \ref{lipschitz_assumption} with $\phi=0$ and $\psi=0$, we can immediately obtain the existence for general $\phi$ and $\psi$. This will simplify our subsequent proof of existence. Notice that the result relies essentially on the Girsanov Theorem of Proposition \ref{girsanov}.

\begin{Lemma}\label{lemma.phipsi}
Define 
$\overline{g}_t(\omega,y,z,u):=g_t(\omega,y,z,u)-\phi_t(\omega).z-\langle\psi_t(\omega),u\rangle_t.$
Then $(Y,Z,U)$ is a solution of the BSDEJ with generator $g$ and terminal condition $\xi$ under $\mathbb P$ if and only if $( Y,Z,U)$ is a solution of the BSDEJ with generator $\overline g$ and terminal condition $\xi$ under $\mathbb Q$ where
$$\frac{d\mathbb Q}{d\mathbb P}=\mathcal E\left(\int_0^T\phi_sdB_s+\int_0^T\int_E\psi_s(x)\widetilde\mu(dx,ds)\right).$$
\end{Lemma}

\proof
We have clearly 
\begin{align*}
Y_t\-\-=\-\-\xi\-+\-\int_t^T\overline g_s(Y_s,Z_s,U_s)ds \--\-\int_t^TZ_s(dB_s-\phi_sds)\--\-\int_t^T\-\int_EU_s(x)\-(\widetilde\mu(dx,ds)\--\-\psi_s(x)\nu_s(dx)ds).
\end{align*}

\vspace{0.3em}
Now, by our BMO assumptions on $\phi$ and $\psi$ and the fact that we assumed that $\psi\geq-1+\delta$, we can apply Proposition \ref{girsanov} and $\mathbb Q$ is well defined. Then by Girsanov Theorem, we know that $dB_s-\phi_sds$ and $\widetilde\mu(dx,ds)-\psi_s(x)\nu_s(dx)ds$ are martingales under $\mathbb Q$. Hence the desired result.
\ep

\begin{Remark}
It is clear that if $g$ satisfies Assumption \ref{lipschitz_assumption}, then $\overline g$ defined above satisfies Assumption \ref{lipschitz_assumption} with $\phi=\psi=0$.
\end{Remark}

Following Lemma \ref{lemma.phipsi} we assume for the time being that $g(0,0,0)=\phi=\psi=0$. Our first result is the following

\begin{Theorem}\label{th.small}
Assume that $$\No{\xi}_\infty\leq\frac{1}{2\sqrt{15}\sqrt{2670}\mu e^{\frac32CT}},$$ where $C$ is the Lipschitz constant of $g$ in $y$, and $\mu$ is the constant appearing in Assumption \ref{lipschitz_assumption}. Then under Assumption \ref{lipschitz_assumption} with $\phi=0$, $\psi=0$ and $g(0,0,0)=0$, there exists a unique solution $(Y,Z,U)\in\mathcal S^\infty\times\mathbb H^2_{\rm{BMO}}\times\mathbb J^2_{\rm{BMO}}\cap L^\infty(\nu)$ of the BSDEJ \reff{eq:bsdej}.
\end{Theorem}

\begin{Remark}
Notice that in the above Theorem, we do not need Assumption \ref{assump:hquad}(iii) to hold. This is linked to the fact that, as discussed in Remark \ref{rem.assumptions}, Assumption \ref{lipschitz_assumption} implies a weak version of Assumption \ref{assump:hquad}(iii), which is sufficient for our purpose here.
\end{Remark}

\proof
We first recall that we have with Assumption \ref{lipschitz_assumption} when $g(0,0,0)=\phi=\psi=0$
\begin{equation}\label{es}
\abs{g_t(y,z,u)}\leq C\abs{y}+\mu\abs{z}^2+\mu\No{u}^2_{L^2(\nu_t)}.
\end{equation}

Consider now the map $\Phi:(y,z,u)\in\mathcal S^\infty\times\mathbb H^2_{\rm{BMO}}\times\mathbb J^2_{\rm{BMO}}\cap L^\infty(\nu)\rightarrow (Y,Z,U)$ defined by
\begin{equation}
Y_t=\xi+\int_t^Tg_s(Y_s,z_s,u_s)ds-\int_t^TZ_sdB_s-\int_t^T\int_EU_s(x)\widetilde\mu(dx,ds).
\label{eq:5}
\end{equation}

The above is nothing more than a BSDEJ with jumps whose generator depends only on $Y$ and is Lipschitz. Besides, since $(z,u)\in\mathbb H^2_{\rm{BMO}}\times\mathbb J^2_{\rm{BMO}}\cap L^\infty(\nu)$, using \reff{inte}, \reff{es} and the energy inequalities \reff{energy} we clearly have
$$\mathbb E\left[\left(\int_0^T\abs{g_s(0,z_s,u_s)}ds\right)^2\right]<+\infty.$$

Hence, the existence of $(Y,Z,U)\in\mathcal S^2\times\mathbb H^2\times\mathbb J^2$ is ensured by the results of Barles, Buckdahn and Pardoux \cite{barles} or Li and Tang \cite{li} for Lipschitz BSDEJs with jumps. Of course, we could have let the generator in \reff{eq:5} depend on $(y_s,z_s,u_s)$ instead. The existence of $(Y,Z,U)$ would then have been a consequence of the predictable martingale representation Theorem. However, the form that we have chosen will simplify some of the following estimates.

\vspace{0.3em}
{\bf Step $1$:} We first show that $(Y,Z,U)\in\mathcal S^\infty\times\mathbb H^2_{\rm{BMO}}\times\mathbb J^2_{\rm{BMO}}\cap L^\infty(\nu)$.

\vspace{0.3em}
Recall that by the Lipschitz hypothesis in $y$, there exists a bounded process $\lambda$ such that 
$$g_s(Y_s,z_s,u_s)=\lambda_sY_s+g_s(0,z_s,u_s).$$

Let us now apply It\^o's formula to $e^{\int_t^s\lambda_udu}Y_s$. We obtain easily from Assumption \ref{lipschitz_assumption}
\begin{align*}
Y_t&=\mathbb E_t\left[e^{\int_t^T\lambda_sds}\xi+\int_t^Te^{\int_t^s\lambda_udu}(\lambda_s Y_s+g_s(0,z_s,u_s))ds-\int_t^T\lambda_se^{\int_t^s\lambda_udu}Y_sds\right]\\
&\leq \mathbb E_t\left[e^{\int_t^T\lambda_sds}\xi+\mu\int_t^Te^{\int_t^s\lambda_udu}\left(\abs{z_s}^2+\int_Eu_s^2(x)\nu_s(dx)\right)ds\right]\\
&\leq \No{\xi}_\infty+C\left(\No{z}^2_{\mathbb H^2_{\rm{BMO}}}+\No{u}^2_{\mathbb J^2_{\rm{BMO}}}\right).
\end{align*}

Therefore $Y$ is bounded and consequently, since its jumps are also bounded, we know that there is a version of $U$ such that 
$$\No{U}_{L^\infty(\nu)}\leq 2\No{Y}_{\mathcal S^\infty}.$$

Let us now prove that $(Z,U)\in \mathbb H^2_{\rm{BMO}}\times\mathbb J^2_{\rm{BMO}}$. Applying It\^o's formula to $e^{\eta t}\abs{Y_t}^2$ for some $\eta>0$, we obtain for any stopping time $\tau\in\mathcal T_0^T$
\begin{align*}
&e^{\eta \tau}\abs{Y_\tau}^2+\mathbb E_\tau\left[\int_\tau^Te^{\eta s}\abs{Z_s}^2ds+\int_\tau^T\int_Ee^{\eta s}U_s^2(x)\nu_s(dx)ds\right]\\
&=\mathbb E_\tau\left[e^{\eta T}\xi^2+2\int_\tau^Te^{\eta s}Y_sg_s(Y_s,z_s,u_s)ds-\eta\int_\tau^Te^{\eta s}\abs{Y_s}^2ds\right]\\
&\leq \mathbb E_\tau\left[e^{\eta T}\xi^2+(2C-\eta)\int_\tau^Te^{\eta s}\abs{Y_s}^2ds+2\No{Y}_{\mathcal S^\infty}\int_\tau^Te^{\eta s}\abs{g_s(0,z_s,u_s)}ds\right].
\end{align*}

Choosing $\eta= 2C$, and using the elementary inequality $2ab\leq \frac{a^2}{\eps}+\eps b^2$, we obtain
\begin{align*}
&\abs{Y_\tau}^2+\mathbb E_\tau\left[\int_\tau^T\abs{Z_s}^2ds+\int_\tau^T\int_EU_s^2(x)\nu_s(dx)ds\right]\\
&\leq \mathbb E_\tau\left[e^{\eta T}\xi^2+\eps\No{Y}_{\mathcal S^\infty}^2+\frac{e^{2\eta T}}{\eps}\left(\int_\tau^T\abs{g_s(0,z_s,u_s)}ds\right)^2\right].
\end{align*}

Hence,
\begin{align*}
&(1-\eps)\No{Y}_{\mathcal S^\infty}^2+\No{Z}^2_{\mathbb H^2_{\rm{BMO}}}+\No{U}^2_{\mathbb J^2_{\rm{BMO}}}\leq e^{\eta T}\No{\xi}_\infty^2+64\mu^2\frac{e^{2\eta T}}{\eps}\left(\No{z}_{\mathbb H^2_{\rm{BMO}}}^4+\No{u}_{\mathbb J^2_{\rm{BMO}}}^4\right).
\end{align*}

And finally, choosing $\eps=1/2$
$$\No{Y}_{\mathcal S^\infty}^2+\No{Z}^2_{\mathbb H^2_{\rm{BMO}}}+\No{U}^2_{\mathbb J^2_{\rm{BMO}}}\leq 2e^{\eta T}\No{\xi}_\infty^2+256\mu^2e^{2\eta T}\left(\No{z}_{\mathbb H^2_{\rm{BMO}}}^4+\No{u}_{\mathbb J^2_{\rm{BMO}}}^4\right).$$

Our problem now is that the norms for $Z$ and $U$ in the left-hand side above are to the power $2$, while they are to the power $4$ on the right-hand side. Therefore, it will clearly be impossible for us to prevent an explosion if we do not first start by restricting ourselves in some ball with a well chosen radius. This is exactly the mathematical manifestation of the phenomenon discussed at the beginning of this section. Define therefore $R=\frac{1}{2\sqrt{2670}\mu e^{\eta T}}$, and assume that $\No{\xi}_\infty\leq \frac{R}{\sqrt{15}e^{\frac12\eta T}}$ and that
$$\No{y}^2_{\mathcal S^\infty}+\No{z}^2_{\mathbb H^2_{\rm{BMO}}}+\No{u}^2_{\mathbb J^2_{\rm{BMO}}}+\No{u}^2_{L^\infty(\nu)}\leq R^2.$$

Denote $\Lambda:=\No{Y}^2_{\mathcal S^\infty}+\No{Z}^2_{\mathbb H^2_{\rm{BMO}}}+\No{U}^2_{\mathbb J^2_{\rm{BMO}}}+\No{U}^2_{L^\infty(\nu)}$. We have, since $\No{U}^2_{L^\infty(\nu)}\leq 4\No{Y}^2_{\mathcal S^\infty}$
\begin{align*}
\Lambda\leq 5\No{Y}^2_{\mathcal S^\infty}+\No{Z}^2_{\mathbb H^2_{\rm{BMO}}}+\No{U}^2_{\mathbb J^2_{\rm{BMO}}}&\leq 10e^{\eta T}\No{\xi}_\infty^2+1280\mu^2e^{2\eta T}\left(\No{z}_{\mathbb H^2_{\rm{BMO}}}^4+\No{u}_{\mathbb J^2_{\rm{BMO}}}^4\right)\\
&\leq \frac{2R^2}{3}+3560\mu^2e^{2\eta T}R^4=\frac{2R^2}{3}+ \frac{R^2}{3}=R^2.
\end{align*}

Hence if $\mathcal B_R$ is the ball of radius $R$ in $\mathcal S^\infty\times\mathbb H^2_{\rm{BMO}}\times\mathbb J^2_{\rm{BMO}}\cap L^\infty(\nu)$, we have shown that 
$\Phi(\mathcal B_R)\subset \mathcal B_R$.

\vspace{0.3em}
{\bf Step $2$:} We show that $\Phi$ is a contraction in this ball of radius $R$.

\vspace{0.3em}
For $i=1,2$ and $(y^i,z^i,u^i)\in\mathcal B_R$, we denote $(Y^i,Z^i,U^i):=\Phi(y^i,z^i,u^i)$ and
\begin{align*}
&\delta y:=y^1-y^2, \ \delta z:=z^1-z^2,\ \delta u:=u^1-u^2,\ \delta Y:=Y^1-Y^2\\
&\delta Z:=Z^1-Z^2,\ \delta U:=U^1-U^2,\ \delta g:=g(Y^2,z^1,u^1)-g(Y^2,z^2,u^2).
\end{align*}

Arguing as above, we obtain easily
$$\No{\delta Y}_{\mathcal S^\infty}^2+\No{\delta Z}^2_{\mathbb H^2_{\rm{BMO}}}+\No{\delta U}^2_{\mathbb J^2_{\rm{BMO}}}\leq 4e^{2\eta T}\underset{\tau\in\mathcal T_0^T}{\sup}\left(\mathbb E_\tau\left[\int_\tau^T\abs{\delta g_s}ds\right]\right)^2.$$

We next estimate that
\begin{align*}
\left(\mathbb E_\tau\left[\int_\tau^T\abs{\delta g_s}ds\right]\right)^2&\leq 2\mu^2\left(\mathbb E_\tau\left[\int_\tau^T\abs{\delta z_s}\left(\abs{z^1_s}+\abs{z^2_s}\right)ds\right]\right)^2\\
&\hspace{0.9em}+2\mu^2\left(\mathbb E_\tau\left[\int_\tau^T\No{\delta u_s}_{L^2(\nu_s)}\left(\No{u_s^1}_{L^2(\nu_s)}+\No{u_s^2}_{L^2(\nu_s)}\right)ds\right]\right)^2\\
&\leq 2\mu^2\left(\mathbb E_\tau\left[\int_\tau^T\abs{\delta z_s}^2ds\right]\mathbb E_\tau\left[\int_\tau^T(\abs{z^1_s}+\abs{z^2_s})^2ds\right]\right.\\
&\hspace{0.9em}\left.+\mathbb E_\tau\left[\int_\tau^T\No{\delta u_s}^2_{L^2(\nu_s)}ds\right]\mathbb E_\tau\left[\int_\tau^T\left(\No{u_s^1}_{L^2(\nu_s)}+\No{u_s^2}_{L^2(\nu_s)}\right)^2ds\right]\right)\\
&\leq 4R^2\mu^2\left(\mathbb E_\tau\left[\int_\tau^T\abs{\delta z_s}^2ds\right]+\mathbb E_\tau\left[\int_\tau^T\int_E\delta u_s^2(x)\nu(dx)ds\right]\right)\\
&\leq 32R^2\mu^2\left(\No{\delta z}_{\mathbb H^2_{\rm{BMO}}}^2+\No{\delta u}_{\mathbb J^2_{\rm{BMO}}}^2\right)
\end{align*}

From these estimates, we obtain, using again the fact that $\No{\delta U}^2_{L^\infty(\nu)}\leq 4\No{\delta Y}_{\mathcal S^\infty}^2$
\begin{align*}
\No{\delta Y}_{\mathcal S^\infty}^2+\No{\delta Z}^2_{\mathbb H^2_{\rm{BMO}}}+\No{\delta U}^2_{\mathbb J^2_{\rm{BMO}}}+\No{\delta U}^2_{L^\infty(\nu)}&\leq 640R^2\mu^2e^{2\eta T}\left(\No{\delta z}_{\mathbb H^2_{\rm{BMO}}}^2+\No{\delta u}_{\mathbb J^2_{\rm{BMO}}}^2\right)\\
&=\frac{16}{267}\left(\No{\delta z}_{\mathbb H^2_{\rm{BMO}}}^2+\No{\delta u}_{\mathbb J^2_{\rm{BMO}}}^2\right).
\end{align*}

Therefore $\Phi$ is a contraction which has a unique fixed point.
\ep

\vspace{0.3em}
Then, from Lemma \ref{lemma.phipsi}, we have immediately the following corollary
\begin{Corollary}\label{corcor}
Assume that $$\No{\xi}_\infty\leq\frac{1}{2\sqrt{15}\sqrt{2670}\mu e^{\frac32CT}},$$ where $C$ is the Lipschitz constant of $g$ in $y$, and $\mu$ is the constant appearing in Assumption \ref{lipschitz_assumption}. Then under Assumption \ref{lipschitz_assumption} with $g(0,0,0)=0$, there exists a unique solution $(Y,Z,U)\in\mathcal S^\infty\times\mathbb H^2_{\rm{BMO}}\times\mathbb J^2_{\rm{BMO}}\cap L^\infty(\nu)$ of the BSDEJ \reff{eq:bsdej}.
\end{Corollary}

We now show how we can get rid of the assumption that $g_t(0,0,0)=0$.
\begin{Corollary}\label{corcor2}
Assume that $$\No{\xi}_\infty+D\No{\int_0^T\abs{g_t(0,0,0)}dt}_\infty\leq\frac{1}{2\sqrt{15}\sqrt{2670}\mu e^{\frac32CT}},$$ where $C$ is the Lipschitz constant of $g$ in $y$, $\mu$ is the constant appearing in Assumption \ref{lipschitz_assumption} and $D$ is a large enough positive constant. Then under Assumption \ref{lipschitz_assumption}, there exists a solution $(Y,Z,U)\in\mathcal S^\infty\times\mathbb H^2_{\rm{BMO}}\times\mathbb J^2_{\rm{BMO}}\cap L^\infty(\nu)$ of the BSDEJ \reff{eq:bsdej}.
\end{Corollary}

\proof
By Corollary \ref{corcor}, we can show the existence of a solution to the BSDEJ with generator $\widetilde g_t(y,z,u):=g_t(y-\int_0^tg_s(0,0,0)ds,z,u)-g_t(0,0,0)$ and terminal condition $\overline \xi:=\xi+\int_0^Tg_t(0,0,0)dt$. Indeed, even though $\widetilde g$ is not null at $(0,0,0)$, it is not difficult to show with the same proof as in Theorem \ref{th.small} that a solution $(\overline Y,\overline Z,\overline U)$ exists (the same type of arguments are used in \cite{tev}). More precisely, $\widetilde g$ still satisfies Assumption \ref{lipschitz_assumption}(i) and when $\phi$ and $\psi$ in Assumption \ref{lipschitz_assumption} are equal to $0$, we have the estimate
$$\abs{\widetilde g_t(y,z,u)}\leq C\No{\int_0^T\abs{g_s(0,0,0)}ds}_\infty +C\abs{y}+\mu\abs{z}^2+\mu\No{u}^2_{L^2(\nu_t)},$$
which is the counterpart of \reff{es}. Thus, since the constant term in the above estimate is assumed to be small enough, it will play the same role as $\No{\xi}_\infty$ in the first Step of the proof of Theorem \ref{th.small}. 

\vspace{0.3em}
For the Step $2$, everything still work thanks to the following estimate
\begin{align*}
\abs{\widetilde g_t(Y^2,z^1,u^1)-\widetilde g_t(Y^2,z^2,u^2)}\leq& \ \mu\abs{z^1-z^2}\left(\abs{z^1}+\abs{z^2}\right)\\
&+\mu\No{u^1-u^2}_{L^2(\nu_t)}\left(\No{u^1}_{L^2(\nu_t)}+\No{u^2}_{L^2(\nu_t)}\right).
\end{align*}

Then, if we define $(Y_t,Z_t,U_t):=(\overline Y_t-\int_0^tg_s(0,0,0)ds,\overline Z_t,\overline U_t),$ it is clear that it is a solution to the BSDEJ with generator $g$ and terminal condition $\xi$.
\ep

\begin{Remark}
We emphasize that the above proof of existence extends readily to a terminal condition which is in $\mathbb R^n$ for any $n\geq 2$.
\end{Remark}

\section{Existence for a bounded terminal condition}\label{sec.ex2}
We now show that we can still prove existence of a solution for any bounded terminal condition. In return, we will now have to strengthen once more our assumptions on the generator. Intuitively speaking, the Lipschitz and local Lipschitz assumptions in Assumption \ref{lipschitz_assumption} are no longer enough and are replaced by stronger regularity assumptions.

\begin{Assumption}\label{assump:hh}
\begin{itemize}
\item[\rm{(i)}] $g$ is uniformly Lipschitz in $y$.
$$\abs{g_t(\omega,y,z,u)-g_t(\omega,y',z,u)}\leq C\abs{y-y'}\text{ for all }(\omega,t,y,y',z,u).$$
\item[\rm{(ii)}] $g$ is $C^2$ in $z$ and there are $\theta>0$ and $(r_t)_{0\leq t\leq T}\in\mathbb H^2_{\rm{BMO}}$,  such that for all $(t,\omega,y,z,u)$,
$$\lvert D_z g_t(\omega,y,z,u)\rvert\leq r_t + \theta\abs{z}, \ \lvert D^2_{zz} g_t(\omega,y,z,u)\rvert\leq \theta.$$
\item[\rm{(iii)}] $g$ is twice Fr\'echet differentiable in the Banach space $L^2(\nu)$ and there are constants $\theta$, $\delta>0$, $C_1\geq-1+\delta$, $C_2\geq 0$ and a predictable function $m\in\mathbb J^2_{\rm{BMO}}$ s.t. for all $(t,\omega,y,z,u,x)$,
$$\abs{D_u g_t(\omega,y,z,u)}\leq m_t+\theta\abs{u},\ C_1(1\wedge\abs{x})\leq D_ug_t(\omega,y,z,u)(x)\leq C_2(1\wedge\abs{x})$$
$$ \No{D^2_u g_t(\omega,y,z,u)}_{L^2(\nu_t)}\leq \theta.$$
\end{itemize}
\end{Assumption}

\begin{Remark}
The assumptions (ii) and (iii) above are generalizations to the jump case of the assumptions considered by Tevzadze \cite{tev}. They will only be useful in our proof of existence and are tailor-made to allow us to apply the Girsanov transformation of Proposition \ref{girsanov}. Notice also that since the space $L^2(\nu)$ is clearly a Banach space, there is no problem to define the Fr\'echet derivative.
\end{Remark}

We emphasize here that Assumption \ref{assump:hh} is stronger than Assumption \ref{lipschitz_assumption}. Indeed, we have the following result

\begin{Lemma}\label{lemma.assump}
If Assumption \ref{assump:hh}(ii) and (iii) hold, then so do Assumption \ref{lipschitz_assumption}(ii) and (iii).
\end{Lemma}

\proof
We will only show that if Assumption \ref{assump:hh}$\rm{(iii)}$ holds, so does Assumption \ref{lipschitz_assumption}$\rm{(iii)}$, the proof being similar for Assumption \ref{assump:hh}$\rm{(ii)}$. Since $g$ is twice Fr\'echet differentiable in $u$, we introduce the process $\psi_t:=D_ug_t(y,z,0)$ which is bounded from above by $m$ and from below by $C_1\geq-1+\delta$ by assumption. Thus, $\psi\in\mathbb J^2_{\rm{BMO}}$. By the mean value theorem, we compute that for some $\lambda\in[0,1]$ and with $u_\lambda:=\lambda u+(1-\lambda)u'$
\begin{align*}
\abs{g_t(y,z,u)-g_t(y,z,u')-\langle\psi_t,u-u'\rangle_t}&\leq\No{D_ug_t(y,z,u_\lambda)-\psi_t}\No{u-u'}_{L^2(\nu_t)}\\
&\leq \theta\No{\lambda u+(1-\lambda)u'}_{L^2(\nu_t)}\No{u-u'}_{L^2(\nu_t)},
\end{align*}
by the bound on $D^2_ug$. The result now follows easily.
\ep

\vspace{0.3em}
We can now state our main existence result.

\begin{Theorem}
Let $\xi\in\mathbb L^\infty$. Under Assumptions \ref{assump:hquad} and \ref{assump:hh}, there exists a solution $(Y,Z,U)\in\mathcal S^\infty\times\mathbb H^2_{\rm{BMO}}\times\mathbb J^2_{\rm{BMO}}\cap L^\infty(\nu)$ of the BSDEJ \reff{eq:bsdej}.
\end{Theorem}

\begin{Remark}
Of course, our existence results are, somehow, less general than the ones obtained in \cite{elmatn} and \cite{ngou}, since they consider generators which only satisfy Assumption \ref{assump:hquad} and are continuous. Moreover, their terminal conditions are not necessarily bounded. However, we emphasize that in the case of a small terminal condition, our result allows us to no longer assume the structure condition of Assumption \ref{assump:hquad}, which can be restrictive from the point of view of applications. Notwithstanding this, we would also like to remind the reader that our approach is fundamentally different from theirs, and allows to obtain solutions from Picard iterations. This property  could be useful for numerical simulations.
\end{Remark}

The idea of the proof is to find a "good" splitting of the BSDEJ into the sum of BSDEJs for which the terminal condition is small and existence holds. Then we paste everything together. This is during this pasting step that the regularity of the generator in $z$ and $u$ in Assumption \ref{assump:hh} is going to be important.

\vspace{0.3em}
\proof

$\rm{(i)}$ We first assume that $g_t(0,0,0)=0$. Consider an arbitrary decomposition of $\xi$
$$\xi=\sum_{i=1}^n\xi_i\text{ such that }\No{\xi_i}_{\infty}\leq \frac{1}{2\sqrt{15}\sqrt{2670}\mu e^{\frac32CT}},\text{ for all $i$.}$$

We will now construct a solution to \reff{eq:bsdej} recursively.

\vspace{0.3em}
{\bf Step $1$} We define $g^1:=g$ and $(Y^1,Z^1,U^1)$ as the unique solution of
\begin{equation}
Y_t^1=\xi_1+\int_t^Tg_s^1(Y_s^1,Z_s^1,U_s^1)ds-\int_t^TZ_s^1dB_s-\int_t^T\int_EU^1_s(x)\widetilde{\mu}(ds,dx),\ \mathbb P-a.s.
\label{eq:bsdej1}
\end{equation}

Let us show why this solution exists. Since $g^1$ satisfies Assumption \ref{assump:hh}, we know by Lemma \ref{lemma.assump} that it satisfies Assumption \ref{lipschitz_assumption} with $\phi_t:=D_zg_t(y,0,u)$ and $\psi_t:=D_ug_t(y,z,0)$, these processes being respectively in $\mathbb H^2_{\rm{BMO}}$ and $\mathbb J^2_{\rm{BMO}}\cap L^\infty(\nu)$ by assumption. Furthermore, we have $\psi_t(x)\geq C_1(1\wedge\abs{x})$ with $C_1\geq -1+\delta$. Thanks to Theorem \ref{th.small} and with the notations of Lemma \ref{lemma.phipsi}, we can then define the solution to the BSDEJ with driver $\overline g^1$ (which still satisfies $\overline g^1(0,0,0)=0$) and terminal condition $ \xi_1$ under the probability measure $\mathbb Q^1$ defined by
$$\frac{d\mathbb Q^1}{d\mathbb P}=\mathcal E\left(\int_0^T\phi_sdB_s+\int_0^T\int_E\psi_s(x)\widetilde\mu(dx,ds)\right).$$

Thanks to Lemma \ref{lemma.phipsi}, this gives us a solution $(Y^1,Z^1,U^1)$ to \reff{eq:bsdej1} with $Y^1$ bounded, which in turn implies with Lemma \ref{lemma.bmo} that $(Z^1,U^1)\in\mathbb H^2_{\rm{BMO}}\times\mathbb J^2_{\rm{BMO}}\cap L^\infty(\nu)$.

\vspace{0.3em}
{\bf Step 2} We assume that we have constructed similarly $(Y^j,Z^j,U^j)\in\mathcal S^\infty\times\mathbb H^2_{\rm{BMO}}\times\mathbb J^2_{\rm{BMO}}\cap L^\infty(\nu)$ for $j\leq i-1$. We then define the generator
$$g_t^i(y,z,u):=g_t\left(\overline{Y}^{i-1}_t+y,\overline{Z}^{i-1}_t+z,\overline{U}^{i-1}_t+u\right)-g_t\left(\overline{Y}^{i-1}_t,\overline{Z}^{i-1}_t,\overline{U}^{i-1}_t\right),$$
where
$$\overline{Y}^{i-1}_t:=\sum_{j=1}^{i-1}Y_t^j,\ \overline{Z}^{i-1}_t:=\sum_{j=1}^{i-1}Z_t^j,\ \overline{U}^{i-1}_t:=\sum_{j=1}^{i-1}U_t^j.$$

Notice that $g^i(0,0,0)=0$ and since $g$ satisfies Assumption \ref{assump:hquad}(iii), we have the estimate
\begin{align*}
&g^i_t(y,z,u)\leq 2\alpha_t+\beta\abs{y+\overline{Y}_t^{i-1}}+\beta\abs{\overline{Y}_t^{i-1}}+\gamma\abs{z+\overline{Z}^{i-1}_t}^2+\gamma\abs{\overline{Z}^{i-1}_t}^2\\
&\hspace{5.1em} +\frac1\gamma j_t\left(\gamma\left(u+\overline{U}^{i-1}_t\right)\right)+\frac1\gamma j_t\left(\gamma\overline{U}^{i-1}_t\right)\\
&\leq2\alpha_t+2\beta\abs{\overline{Y}_t^{i-1}}+3\gamma\abs{\overline{Z}^{i-1}_t}^2+\frac{1}{\gamma} j_t\left(\gamma\overline{U}^{i-1}_t\right)+\frac{1}{2\gamma} j_t\left(2\gamma\overline{U}^{i-1}_t\right)+\beta\abs{y}+2\gamma\abs{z}^2+\frac{1}{2\gamma} j_t\left(2\gamma u\right),
\end{align*}
where we used the inequalities $(a+b)^2\leq 2(a^2+b^2)$ and the fact that for all $(\gamma_1,\gamma_2)\in\mathbb R^2$ (see \cite{kpz4} for a proof)
\begin{equation}\label{gammma}
(\gamma_1+\gamma_2)\left(e^{\frac{x+y}{\gamma_1+\gamma_2}}-1- \frac{x+y}{\gamma_1+\gamma_2}\right)\leq \gamma_1\left(e^{\frac{x}{\gamma_1}}-1- \frac{x}{\gamma_1}\right) + \gamma_2\left(e^{\frac{y}{\gamma_2}}-1- \frac{y}{\gamma_2}\right).
\end{equation}

\vspace{0.3em}
Then, since $(\overline{Y}^{i-1},\overline{Z}^{i-1},\overline{U}^{i-1})\in\mathcal S^\infty\times\mathbb H^2_{\rm{BMO}}\times\mathbb J^2_{\rm{BMO}}\cap L^\infty(\nu)$ , we know that the term which does not depend on $(y,z,u)$ above satisfies the same integrability condition as $g_t(0,0,0)$ in \reff{inte} (see also the arguments we used in Remark \ref{rem.assumptions}). Therefore, since $g^i(0,0,0)=0$, we have one side of the inequality in Assumption \ref{assump:hquad}(iii), and the other one can be proved similarly. This yields that $g^i$ satisfies Assumption \ref{assump:hquad}.

\vspace{0.3em}
Similarly as in Step $1$, we will now show that there exists a solution $(Y^i,Z^i,U^i)\in\mathcal S^\infty\times\mathbb H^2_{\rm{BMO}}\times\mathbb J^2_{\rm{BMO}}\cap L^\infty(\nu)$ to the BSDEJ
\begin{equation}
Y_t^i=\xi_i+\int_t^Tg_s^i(Y_s^i,Z_s^i,U_s^i)ds-\int_t^TZ_s^idB_s-\int_t^T\int_EU_s^i(x)\widetilde\mu(dx,ds),\ \mathbb P-a.s.
\label{eq:bsdeji}
\end{equation}

Since $g$ satisfies Assumptions \ref{assump:hh}, we can define
\begin{align*}
\phi^i_t:=D_zg^i_t(y,0,u)=D_zg_t(\overline{Y}^{i-1}_t+y,\overline{Z}^{i-1}_t,\overline{U}^{i-1}_t+u),\\ 
\psi^i_t:=D_ug^i_t(y,z,0)=D_ug_t(\overline{Y}^{i-1}_t+y,\overline{Z}^{i-1}_t+z,\overline{U}^{i-1}_t).
\end{align*}

We then know that
$$\abs{\phi^i_t}\leq r_t+\theta\abs{\overline{Z}^{i-1}_t},\ \abs{\psi^i_t}\leq m_t+\theta\abs{\overline{U}^{i-1}_t}, \ \psi_t^i(x)\geq C_1(1\wedge\abs{x})\geq-1+\delta.$$

Since by hypothesis $(\overline{Z}^{i-1}, \overline{U}^{i-1})\in\mathbb H^2_{\rm{BMO}}\times\mathbb J^2_{\rm{BMO}}\cap L^\infty(\nu)$, we can define $\mathbb Q^i$ by
$$\frac{d\mathbb Q^i}{d\mathbb P}=\mathcal E\left(\int_0^T\phi_s^idB_s+\int_0^T\int_E\psi_s^i(x)\widetilde\mu(dx,ds)\right).$$

Now, using the notations of Lemma \ref{lemma.phipsi}, we define a generator $\bar g^i$ from $g^i$ (which still satisfies $\overline g^i(0,0,0)=0$). It is then easy to check that $\bar g^i$ satisfies Assumption \ref{lipschitz_assumption}. Therefore, by Theorem \ref{th.small}, we obtain the existence of a solution to the BSDEJ with generator $\bar g^i$ and terminal condition $\xi_i$ under $\mathbb Q^i$. Using Lemma \ref{lemma.phipsi}, this provides a solution $(Y^i,Z^i,U^i)$ with $Y^i$ bounded to the BSDEJ \reff{eq:bsdeji}. By Lemma \ref{lemma.bmo} and since $g^i$ satisfies Assumption \ref{assump:hquad}, the boundedness of $Y^i$ implies that $(Z^i,U^i)\in\mathbb H^2_{\rm{BMO}}\times\mathbb J^2_{\rm{BMO}}\cap L^\infty(\nu)$ and therefore that $(\overline{Y}^i,\overline{Z}^i,\overline{U}^i)\in\mathcal S^\infty\times\mathbb H^2_{\rm{BMO}}\times\mathbb J^2_{\rm{BMO}}\cap L^\infty(\nu)$.

\vspace{0.3em}
{\bf Step $3$} Finally, by summing the BSDEJs \reff{eq:bsdeji}, we obtain
$$\overline{Y}^n=\xi+\int_t^Tg_s(\overline{Y}^n_s,\overline{Z}^n_s,\overline{U}^n_s)ds-\int_t^T\overline{Z}^n_sdB_s-\int_t^T\int_E\overline{U}^n_s(x)\widetilde\mu(dx,ds).$$

Since $\overline{Y}^n$ is bounded (because the $Y^i$ are all bounded), Lemma \ref{lemma.bmo} implies that $(\overline{Z}^n,\overline{U}^n)\in\mathbb H^2_{\rm{BMO}}\times\mathbb J^2_{\rm{BMO}}\cap L^\infty(\nu)$, which ends the proof.

\vspace{0.3em}
$\rm{(ii)}$ In the general case $g_t(0,0,0)\neq 0$, we can argue exactly as in Corollary \ref{corcor2} (see also Proposition $2$ in \cite{tev}) to obtain the result.
\ep

\section{Comparison and stability}\label{sec.comp}
\subsection{A uniqueness result}
We emphasize that the above theorems provide an existence result for every bounded terminal condition, but we only have uniqueness when the infinite norm of $\xi$ is small enough. In order to have a general uniqueness result, we add the following assumptions, which were first introduced by Royer \cite{roy} and Briand and Hu \cite{bh2}. Notice that \cite{ngou} also considers Assumption \ref{assump.roy} to recover uniqueness.
\begin{Assumption}\label{assump.roy}
For every $(y,z,u,u')$ there exists a predictable and $\mathcal E$-measurable process $(\gamma_t)$ such that
$$g_t(y,z,u)-g_t(y,z,u')\leq\int_E\gamma_t(x)(u-u')(x)\nu_t(dx),$$
where there exist constants $C_2>0$ and $C_1\geq-1+\delta$ for some $\delta>0$ such that
$$C_1(1\wedge\abs{x})\leq\gamma_t(x)\leq C_2(1\wedge\abs{x}).$$
\end{Assumption}

\begin{Assumption}\label{assump.bh}
$g$ is jointly convex in $(z,u)$.
\end{Assumption}
We then have the following result
\begin{Theorem}\label{th.unique}
Assume that $\xi\in\mathbb L^\infty$, and that the generator $g$ satisfies either 
\begin{enumerate}
	\item[\rm{(i)}] Assumptions \ref{assump:hquad}, \ref{assump:hh}(i),(ii) and \ref{assump.roy}.
	\item[\rm{(ii)}] Assumptions \ref{assump:hquad}, \ref{assump:hh} and \ref{assump.bh}, and that $g(0,0,0)$ and the process $\alpha$ appearing in Assumption \ref{assump:hquad}(iii) are bounded by some constant $M>0$.
\end{enumerate}
Then there exists a unique solution to the BSDEJ \reff{eq:bsdej}.
\end{Theorem}

In order to prove this Theorem, we will use the following comparison Theorem for BSDEJs 

\begin{Proposition}\label{prop.comp}
Let $\xi^1$ and $\xi^2$ be two $\mathcal F_T$-measurable random variables. Let $g^1$ be a function satisfying either of the following
\begin{enumerate}
	\item[\rm{(i)}] Assumptions \ref{assump:hquad}, \ref{lipschitz_assumption}(i),(ii) and \ref{assump.roy}.
	\item[\rm{(ii)}] Assumptions \ref{assump:hquad}, \ref{lipschitz_assumption}(i)  and \ref{assump.bh}, and that $\abs{g^1(0,0,0)}+\alpha\leq M$ where $\alpha$ is the process appearing in Assumption \ref{assump:hquad}(iii) and $M$ is a positive constant.
\end{enumerate}

Let $g^2$ be another function and for $i=1,2$, let $(Y^i,Z^i,U^i)$ be the solution of the BSDEJ with terminal condition $\xi^i$ and generator $g^i$ (we assume that existence holds in our spaces), that is to say for every $t\in[0,T]$
$$Y^i_t=\xi^i+\int_t^Tg^i_s(Y^i_s,Z^i_s,U^i_s)ds-\int_t^TZ^i_sdBs-\int_t^T\int_EU^i_s(x)\widetilde\mu(dx,ds),\ \mathbb P-a.s.$$
Assume further that $\xi^1\leq\xi^2,\ \mathbb P-a.s.$ and $g_t^1(Y_t^2,Z_t^2,U_t^2)\leq g_t^2(Y_t^2,Z_t^2,U_t^2),\ \mathbb P-a.s.$ Then $Y_t^1\leq Y_t^2$, $\mathbb P-a.s.$ Moreover in case (i), if in addition we have $Y^1_0=Y^2_0$, then for all $t$, $Y^1_t=Y^2_t$, $Z_t^1=Z_t^2$ and $U_t^1=U_t^2$, $\mathbb P-a.s.$
\end{Proposition}

\begin{Remark}
Of course, we can replace the convexity property in Assumption \ref{assump.bh} by concavity without changing the results of Proposition \ref{prop.comp}. Indeed, if  $Y$ is a solution to the BSDEJ with convex generator $g$ and terminal condition $\xi$, then $-Y$ is a solution to the BSDEJ with concave generator $\widetilde g(y,z,u):=-g(-y,-z,-u)$ and terminal condition $-\xi$. then we can apply the results of Proposition \ref{prop.comp}. 
\end{Remark}

\proof
{\bf Step $1$} In order to prove $(i)$, let us note 
\begin{align*}
&\delta Y:=Y^1-Y^2,\ \delta Z:=Z^1-Z^2,\ \delta U:=U^1-U^2,\ \delta \xi:=\xi^1-\xi^2\\
&\delta g_t:=g_t^1(Y_t^2,Z_t^2,U_t^2)- g_t^2(Y_t^2,Z_t^2,U_t^2).
\end{align*}

Using Assumption \ref{lipschitz_assumption}$\rm{(i)},\rm{(ii)}$, we know that there exist a bounded $\lambda$ and a process $\eta$ with
\begin{equation}
\abs{\eta_s}\leq\mu\left(\abs{Z_s^1}+\abs{Z_s^2}\right),
\label{eq:machin}
\end{equation}
such that
\begin{align}\label{ew}
\nonumber\delta Y_t&=\delta\xi+\int_t^T\delta g_sds+\int_t^T\lambda_s\delta Y_sds+\int_t^T(\eta_s+\phi_s)\delta Z_sds\\
\nonumber&\hspace{0.9em}+\int_t^Tg_s^1(Y_s^1,Z_s^1,U_s^1)-g_s^1(Y_s^1,Z_s^1,U_s^2)ds-\int_t^T\int_E\delta U_s(x)\gamma_s(x)\nu_s(dx)ds\\
&\hspace{0.9em}+\int_t^T\int_E\delta U_s(x)\gamma_s(x)\nu_s(dx)ds-\int_t^T\delta Z_sdB_s-\int_t^T\int_E\delta U_s(x)\widetilde\mu(dx,ds),
\end{align}
where $\gamma$ is the predictable process appearing in the right hand side of Assumption \ref{assump.roy}. 

\vspace{0.3em}
Define for $s\geq t$, $e^{\Lambda_s}:=e^{\int_t^s\lambda_udu},$
and
$$\frac{d\mathbb Q}{d\mathbb P}:=\mathcal E\left(\int_t^s(\eta_s+\phi_s)dB_s+\int_t^s\int_E\gamma_s(x)\widetilde\mu(dx,ds)\right).$$

Since the $Z^i$ are in $\mathbb H^2_{\rm{BMO}}$, so is $\eta$ and by our assumption on $\gamma_s$ the above stochastic exponential defines a true strictly positive uniformly integrable martingale (see Kazamaki \cite{kaz2}). Then applying It\^o's formula and taking conditional expectation under the probability measure $\mathbb Q$, we obtain
\begin{align}
\label{eq.eq}
\nonumber\delta Y_t&=\mathbb E_t^\mathbb Q\left[e^{\Lambda_T} \delta\xi+\int_t^T e^{\Lambda_s}\delta g_sds\right]\\
&\hspace{0.9em}+\mathbb E_t\left[\int_t^T e^{\Lambda_s}\left(g_s^1(Y_s^1,Z_s^1,U_s^1)-g_s^1(Y_s^1,Z_s^1,U_s^2)-\int_E\gamma_s(x)\delta U_s(x)\nu_s(dx)\right)ds\right]\leq 0,
\end{align}
using Assumption \ref{assump.roy}.

\vspace{0.3em}
{\bf Step $2$} The proof of the comparison result when (ii) holds is a generalization of Theorem $5$ in \cite{bh2}. However, due to the presence of jumps our proof is slightly different. For the convenience of the reader, we will highlight the main differences during the proof.

\vspace{0.3em}
For any $\theta\in(0,1)$ let us denote
$$\delta Y_t:=Y_t^1-\theta Y_t^2,\ \delta Z_t:=Z_t^1-\theta Z_t^2,\ \delta U_t:=U_t^1-\theta U_t^2,\ \delta\xi:=\xi^1-\theta\xi^2.$$

First of all, we have for all $t\in[0,T]$
$$\delta Y_t=\delta\xi+\int_t^TG_sds-\int_t^T\delta Z_sdB_s-\int_t^T\int_E\delta U_s(x)\widetilde\mu(dx,ds),\ \mathbb P-a.s.,$$
where
$$G_t:=g_t^1(Y_t^1,Z_t^1,U_t^1)-\theta g_t^2(Y_t^2,Z_t^2,U_t^2).$$

We emphasize that unlike in \cite{bh2}, we have not linearized the generator in $y$ using the Assumption \ref{lipschitz_assumption}$\rm{(i)}$. It will be clear later on why.

\vspace{0.3em}
We will now bound $G_t$ from above. First, we rewrite it as
$$G_t=G_t^1+G_t^2+G_t^3,$$
where
\begin{align*}
&G_t^1:=g_t^1(Y_t^1,Z_t^1,U_t^1)-g_t^1(Y_t^2,Z_t^1,U_t^1), \ G_t^2:=g_t^1(Y_t^2,Z_t^1,U_t^1)-\theta g_t^1(Y_t^2,Z_t^2,U_t^2)\\ &G_t^3:=\theta\left(g_t^1(Y_t^2,Z_t^2,U_t^2)-g_t^2(Y_t^2,Z_t^2,U_t^2)\right).
\end{align*}

Then, we have using Assumption \ref{lipschitz_assumption}(i)
\begin{align}\label{G1}
\nonumber G_t^1&=g_t^1(Y_t^1,Z_t^1,U_t^1)-g_t^1(\theta Y_t^2,Z_t^1,U_t^1)+g_t^1(\theta Y_t^2,Z_t^1,U_t^1)-g_t^1(Y_t^2,Z_t^1,U_t^1)\\
&\leq C\left(\abs{\delta y_t}+(1-\theta)\abs{y_t^2}\right).
\end{align}

Next, we estimate $G^2$ using Assumption \ref{assump:hquad} and the convexity in $(z,u)$ of $g^1$ \begin{align*}
\nonumber g_t^1(Y_t^2,Z_t^1,U_t^1)&=g_t^1\left(Y_t^2,\theta Z_t^2+(1-\theta)\frac{Z_t^1-\theta Z_t^2}{1-\theta},\theta U_t^2+(1-\theta)\frac{U_t^1-\theta U_t^2}{1-\theta}\right)\\
\nonumber &\leq \theta g_t^1(Y_t^2,Z_t^2,U_t^2)+(1-\theta)g_t^1\left(Y_t^2,\frac{\delta Z_t}{1-\theta},\frac{\delta U_t}{1-\theta}\right)\\
&\leq \theta g_t^1(Y_t^2,Z_t^2,U_t^2)+(1-\theta)\left(M+\beta\abs{Y_t^2}\right)+\frac{\gamma}{2(1-\theta)}\abs{\delta Z_t}^2\\
&\hspace{0.9em}+\frac{1-\theta}{\gamma}j_t\left(\frac{\gamma}{1-\theta}\delta U_t\right).
\end{align*}

Hence
\begin{equation}
G_t^2\leq (1-\theta)\left(M+\beta\abs{Y_t^2}\right)+\frac{\gamma}{2(1-\theta)}\abs{\delta Z_t}^2+\frac{1-\theta}{\gamma}j_t\left(\frac{\gamma}{1-\theta}\delta U_t\right).
\label{eq:G2}
\end{equation}

Finally, $G^3$ is negative by assumption. Therefore, using \reff{G1} and \reff{eq:G2}, we obtain
\begin{equation}
G_t\leq C\abs{\delta Y_t}+(1-\theta)\left(M+\widetilde\beta\abs{Y_t^2}\right)+\frac{\gamma}{2(1-\theta)}\abs{Z_t}^2+\frac{1-\theta}{\gamma}j_t\left(\frac{\gamma}{1-\theta}\delta U_t\right),
\label{eq:G3}
\end{equation}
where $\widetilde\beta:=\beta +C.$

\vspace{0.3em}
Now we will get rid of the quadratic and exponential terms in $z$ and $u$ using a classical exponential change. Let us then denote for some $\nu>0$
\begin{align*}
&P_t:=e^{\nu\delta Y_t},\ Q_t:=\nu e^{\nu \delta Y_t}\delta Z_t, \ R_t(x):=e^{\nu \delta Y_{t^-}}\left(e^{\nu\delta U_t(x)}-1\right).
\end{align*}

\vspace{0.3em}
By It\^o's formula we obtain for every $t\in[0,T]$, $\mathbb P-a.s.$
$$P_t=P_T+\int_t^T\nu P_s\left(G_s-\frac\nu2\abs{\delta Z_s}^2-\frac1\nu j_s(\nu\delta U_s)\right)ds-\int_t^TQ_sdB_s-\int_t^T\int_ER_s(x)\widetilde\mu(dx,ds).$$

Now choose $\nu=\gamma/(1-\theta)$. We emphasize that this is here that the presence of jumps forces us to change our proof in comparison with the one in \cite{bh2}. Indeed, if we had immediately linearized in $y$ then we could not have chosen $\nu$ constant such that the quadratic and exponentials terms in \reff{eq:G3} would disappear. This is not a problem in \cite{bh2}, since they can choose $\nu$ of the form $M/(1-\theta)$ with $M$ large enough and still make the quadratic term in $z$ disappear. However, in the jump case, the application $\gamma \mapsto \gamma^{-1}j_t(\gamma u)$ is not always increasing, and this trick does not work. Nonetheless, we now define the strictly positive and continuous process
$$D_t:=\exp\left(\gamma\int_0^t\left(M+\widetilde\beta\abs{Y_t^2}+\frac{C}{1-\theta}\abs{\delta Y_s}\right)ds\right).$$

Applying It\^o's formula to $D_tP_t$, we obtain
\begin{align*}
d(D_sP_s)=&-\nu D_sP_s\left(G_s-\frac\nu2\abs{\delta Z_s}^2-\frac{j_s(\nu\delta U_s)}{\nu} -C\abs{\delta Y_s}+(1-\theta)\left(M+\widetilde\beta\abs{Y_s^2}\right)\right)ds\\
&+D_sQ_sdB_s+\int_ED_{s^-}R_s(x)\widetilde\mu(dx,ds).
\end{align*}

Hence, using the inequality \reff{eq:G3}, we deduce
$$D_tP_t\leq \mathbb E_t\left[D_TP_T\right], \ \mathbb P-a.s.,$$
which can be rewritten
$$\delta Y_t\leq \frac{1-\theta}{\gamma}\ln\left(\mathbb E_t\left[\exp\left(\gamma\int_t^T\left(M+\widetilde\beta\abs{Y_t^2}+\frac{C}{1-\theta}\abs{\delta Y_s}\right)ds+\frac{\gamma}{1-\theta}\delta\xi\right)\right]\right),\ \mathbb P-a.s.$$

Next, we have
$$\delta\xi=(1-\theta)\xi^1+\theta\left(\xi^1-\xi^2\right)\leq (1-\theta)\abs{\xi^1}.$$

Consequently, we have for some constant $C_0>0$, independent of $\theta$, using the fact that $Y^2$ and $\xi^1$ are bounded $\mathbb P-a.s.$
\begin{align}
\delta Y_t\leq& \frac{1-\theta}{\gamma}\left(\ln(C_0)+\ln\left(\mathbb E_t\left[\exp\left(\frac{C}{1-\theta}\int_t^T\abs{\delta Y_s}ds\right)\right]\right)\right),\ \mathbb P-a.s.
\label{eq:G5}
\end{align}

We finally argue by contradiction. More precisely, let 
$$\mathcal A:=\left\{\omega\in \Omega, Y_t^1(\omega)>Y_t^2(\omega)\right\},$$
and assume that $\mathbb P(\mathcal A)>0.$ Let us then call $\mathcal N$ the $\mathbb P$-negligible set outside of which \reff{eq:G5} holds. Since $\mathcal A$ has a strictly positive probability, $\mathcal B:=\mathcal A\cap \left(\Omega\backslash\mathcal N\right)$ is not empty and also has a strictly positive probability. Then, we would have from \reff{eq:G5} that for every $\omega \in \mathcal B$ 
\begin{equation}
\delta Y_t(\omega)\leq \frac{1-\theta}{\gamma}\ln(C_0)+\frac{C}{\gamma}\int_t^T\No{\delta Y_s}_{\infty,\mathcal B}ds,
\label{eq:G6}
\end{equation}
where $\No{\cdot}_{\infty,\mathcal B}$ is the usual infinite norm restricted to $\mathcal B$.

\vspace{0.3em}
Now, using the dominated convergence theorem, we can let $\theta\uparrow 1^-$ in \reff{eq:G6} to obtain that for any $\omega \in\mathcal B$
$$ Y_t^1(\omega)-Y_t^2(\omega)\leq \frac{C}{\gamma}\int_t^T\No{ Y_s^1-Y_s^2}_{\infty,\mathcal B}ds,$$
which in turns implies, since $\mathcal B\subset\mathcal A$
$$\No{Y_t^1-Y_t^2}_{\infty,\mathcal B}\leq \frac{C}{\gamma}\int_t^T\No{ Y_s^1-Y_s^2}_{\infty,\mathcal B}ds.$$

But with Gronwall's lemma this implies that $\No{Y_t^1-Y_t^2}_{\infty,\mathcal B}=0$ and the desired contradiction. Hence the result.

\vspace{0.3em}
{\bf Step $3$} Let us now assume that $Y_0^1=Y_0^2$ and that we are in the same framework as in Step $1$. Using this in \reff{eq.eq} above when $t=0$, we obtain 
\begin{align}
\label{eq.eq2}
\nonumber0&=\mathbb E^\mathbb Q\left[e^{\Lambda_T}\delta\xi+\int_0^T e^{\Lambda_s}\delta g_sds\right]\\
&\hspace{0.9em}+\mathbb E^\mathbb Q\left[\int_0^T e^{\Lambda_s}\left(g_s^1(Y_s^1,Z_s^1,U_s^1)-g_s^1(Y_s^1,Z_s^1,U_s^2)-\int_E\delta U_s(x)\gamma_s(x)\nu_s(dx)\right)ds\right]\leq 0.
\end{align}

Hence, since all the above quantities have the same sign, this implies in particular that
$$e^{\Lambda_T}\delta\xi+\int_0^T e^{\Lambda_s}\delta g_sds=0,\ \mathbb P-a.s.$$

Moreover, we also have $\mathbb P-a.s.$
$$\int_0^T e^{\Lambda_s}\left(g_s^1(Y_s^1,Z_s^1,U_s^1)-g_s^1(Y_s^1,Z_s^1,U_s^2)\right)ds=\int_0^T e^{\Lambda_s}\left(\int_E\delta U_s(x)\gamma_s(x)\nu_s(dx)\right)ds.$$

Using this result in \reff{ew}, we obtain with It\^o's formula
\begin{align}
\label{eq.eq3}
\nonumber\delta Y_t&=\int_0^T e^{\Lambda_s}\left(\int_E\delta U_s(x)\gamma_s(x)\nu_s(dx)\right)ds-\int_t^T e^{\Lambda_s}\delta Z_s(dB_s-(\eta_s+\phi_s)ds)\\
&\hspace{0.9em}-\int_t^T\int_E e^{\Lambda_{s^-}}\delta U_s(x)\widetilde\mu(dx,ds).
\end{align}

The right-hand side is a martingale under $\mathbb Q$ with null expectation. Thus, since $\delta Y_t\leq 0$, this implies that $Y_t^1=Y_t^2,$ $\mathbb P-a.s.$ Using this in \reff{eq.eq3}, we obtain that the martingale part must be equal to $0$, which implies that $\delta Z_t=0$ and $\delta U_t=0$. 
\ep

\begin{Remark}\label{rem.imp}
In the above proof of the comparison theorem in case (i), we emphasize that it is actually sufficient that, instead of Assumption \ref{assump.roy}, the generator $g$ satisfies
$$g_s^1(Y_s^1,Z_s^1,U_s^1)-g_s^1(Y_s^1,Z_s^1,U_s^2)\leq\int_E\gamma_s(x)\delta U_s(x)\nu_s(dx),$$
for some $\gamma_s$ such that
$$C_1(1\wedge\abs{x})\leq\gamma_s(x)\leq C_2(1\wedge\abs{x}).$$

Besides, this also holds true for the comparison Theorem for Lipschitz BSDEJs with jumps proved by Royer (see Theorem $2.5$ in \cite{roy}).
\end{Remark}


We can now prove Theorem \ref{th.unique}

\vspace{0.3em}
\proof[Proof of Theorem \ref{th.unique}]
First let us deal with the question of existence.
\begin{enumerate}
	\item[\rm{(i)}] If $g$ satisfies Assumptions \ref{assump:hquad}, \ref{assump:hh}(i),(ii) and \ref{assump.roy}, the existence part can be obtained exactly as in the previous proof, starting from a small terminal condition, and using the fact that Assumption \ref{assump.roy} implies that $g$ is Lipschitz in $u$. Thus we omit it.
	\item[\rm{(ii)}] If $g$ satisfies Assumptions \ref{assump:hquad}, \ref{assump:hh} and \ref{assump.bh}, then we already proved existence for bounded terminal conditions.
\end{enumerate}
The uniqueness is then a simple consequence of the above comparison theorem.
\ep

\begin{Remark}
In \cite{kpz4}, we prove a non-linear Doob-Meyer decomposition and obtain as a consequence a reverse comparison Theorem.
\end{Remark}

\subsection{A priori estimates and stability}
In this subsection, we show that under our hypotheses, we can obtain {\it {\it a priori}} estimates for quadratic BSDEs with jumps. We have the following results

\begin{Proposition}\label{apriori}
Let $(\xi^1,\xi^2)\in\mathbb L^\infty\times\mathbb L^\infty$ and let $g$ be a function satisfying Assumptions \ref{assump:hquad}, \ref{lipschitz_assumption}(i),(ii) and \ref{assump.roy}. Let us consider for $i=1,2$ the solutions $(Y^i,Z^i,U^i)\in\mathcal S^\infty\times\mathbb H^2_{\rm{BMO}}\times\mathbb J^2_{\rm{BMO}}$ of the BSDEJs with generator $g$ and terminal condition $\xi^i$ (once again existence is assumed). Then we have for some constant $C>0$
\begin{align*}
&\No{Y^1-Y^2}_{\mathcal S_\infty}+\No{U^1-U^2}_{L^\infty(\nu)} \leq C\No{\xi^1-\xi^2}_\infty\\
&\No{Z^1-Z^2}^2_{\mathbb H^2_{\rm{BMO}}}+\No{U^1-U^2}^2_{\mathbb J^2_{\rm{BMO}}}\leq C\No{\xi^1-\xi^2}_{\infty}.
\end{align*}
\end{Proposition}

\vspace{0.3em}
\proof
Following exactly the same arguments as in Step $1$ of the proof Proposition \ref{prop.comp}, we obtain with the same notations
$$Y_t^1-Y_t^2=\mathbb E^\mathbb Q_t\left[e^{\Lambda_T}(\xi^1-\xi^2)\right]\leq C\No{\xi^1-\xi^2}_\infty,\ \mathbb P-a.s.$$
Notice then that this implies as usual that there is a version of $(U^1-U^2)$ (still denoted $(U^1-U^2)$ for simplicity) which is bounded by $2\No{Y^1-Y^2}_{\mathcal S_\infty}$. This gives easily the first estimate.

\vspace{0.3em}
Let now $\tau\in\mathcal T_0^T$ be a stopping time. Denote also
$$\delta g_s:=g_s(Y_s^1,Z_s^1,U_s^1)-g_s(Y_s^2,Z_s^2,U_s^2).$$

By It\^o's formula, we have using standard calculations
\begin{align}\label{eq:G7}
\nonumber\mathbb E_\tau\left[\int_\tau^T\abs{Z_s}^2ds+\int_\tau^T\No{U_s}^2_{L^2(\nu_s)}ds\right]\leq&\ \mathbb E_\tau\left[\abs{\xi^1-\xi^2}^2+2\int_\tau^T(Y_s^1-Y_s^2)\delta g_sds\right]\\
\leq& \No{\xi^1-\xi^2}_\infty^2+2\No{Y^1-Y^2}_{\mathcal S_\infty}\mathbb E_\tau\left[\int_\tau^T\abs{\delta g_s}ds\right]. 
\end{align}

Then, using Assumption \ref{assump:hquad}, we estimate
\begin{align*}
\abs{\delta g_t}&\leq C\left(\abs{g_t(0,0,0)}+\alpha_t+\sum_{i=1,2}\abs{Y_t^i}+\abs{Z_t^i}^2+j_t\left(\gamma U^i\right)+j_t\left(-\gamma U^i\right)\right)\\
&\leq C\left(\abs{g_t(0,0,0)}+\alpha_t+\sum_{i=1,2}\abs{Y_t^i}+\abs{Z_t^i}^2+\No{U^i_t}_{L^2(\nu)}^2\right),
\end{align*}
where we used the fact that for every $x$ in a compact subset of $\mathbb R$, $0\leq e^x-1-x\leq Cx^2$. Using this estimate and the integrability assumed on $g_t(0,0,0)$ and $\alpha_t$ in \reff{eq:G7} entails
\begin{align*}
\nonumber&\mathbb E_\tau\left[\int_\tau^T\abs{Z_s}^2ds+\int_\tau^T\No{U_s}^2_{L^2(\nu_s)}ds\right]\\
&\leq \No{\xi^1-\xi^2}_\infty^2+C\No{\xi^1-\xi^2}_{\infty}\left(1+\sum_{i=1,2}\No{Y^i}_{\mathcal S^\infty}+\No{Z^i}_{\mathbb H^2_{\rm{BMO}}}+\No{U^i}_{\mathbb J^2_{\rm{BMO}}}\right)\leq C\No{\xi^1-\xi^2}_\infty, 
\end{align*}
which ends the proof.
\ep

\vspace{0.3em}
\begin{Proposition}\label{apriori2}
Let $(\xi^1,\xi^2)\in\mathbb L^\infty\times\mathbb L^\infty$ and let $g$ be a function satisfying Assumptions \ref{assump:hquad}, \ref{lipschitz_assumption}(i) and \ref{assump.bh} and such that $\abs{g(0,0,0)}+\alpha\leq M$ where $\alpha$ is the process appearing in Assumption \ref{assump:hquad}(iii)  and $M$ is a positive constant. Let us consider for $i=1,2$ the solutions $(Y^i,Z^i,U^i)\in\mathcal S^\infty\times\mathbb H^2_{\rm{BMO}}\times\mathbb J^2_{\rm{BMO}}$ of the BSDEJs with generator $g$ and terminal condition $\xi^i$ (once again existence is assumed). Then we have for some constant $C>0$
\begin{align*}
&\No{Y^1-Y^2}_{\mathcal S^\infty}+\No{U^1-U^2}_{L^\infty(\nu)} \leq C\No{\xi^1-\xi^2}_\infty\\
&\No{Z^1-Z^2}^2_{\mathbb H^2_{\rm{BMO}}}+\No{U^1-U^2}^2_{\mathbb J^2_{\rm{BMO}}}\leq C\No{\xi^1-\xi^2}_{\infty}.
\end{align*}
\end{Proposition}

\vspace{0.3em}
\proof
Following Step $2$ of the proof of Proposition \ref{prop.comp}, we obtain for any $\theta\in(0,1)$
\begin{equation*}
\frac{Y_t^1-\theta Y_t^2}{1-\theta}\leq\frac{1}{\gamma}\ln\left(\mathbb E_t\left[\exp\left(\gamma\int_t^T\left(M+\widetilde\beta\abs{Y_s^2}+\frac{C\abs{Y_s^1-\theta Y_s^2}}{1-\theta}\right)ds+\frac{\gamma(\xi^1-\theta\xi^2)}{1-\theta}\right)\right]\right),
\end{equation*}
and of course by symmetry, the same holds if we interchange the roles of the exponents $1$ and $2$. Since all the quantities above are bounded, we obtain easily after some calculations, after letting $\theta \uparrow 1^-$ and by symmetry
$$\abs{Y_t^1-Y_t^2}\leq C\left(\No{\xi^1-\xi^2}_\infty+\int_t^T\No{Y_s^1-Y_s^2}_\infty ds\right),\ \mathbb P-a.s.$$
%
Hence, we can use Gronwall's lemma to obtain $\No{Y^1-Y^2}_{\mathcal S^\infty}\leq C\No{\xi^1-\xi^2}_\infty.$ All the other estimates can then be obtained as in the proof of Proposition \ref{apriori}.
\ep

\vspace{0.3em}
\paragraph{Acknowledgements:} The authors would like to thank Nicole El Karoui and Anis Matoussi for their precious advices which greatly helped to improve a previous version of this paper.

\end{document}